\documentclass{easychair}

\usepackage{doc}

\usepackage{bbold}

\RequirePackage{scalerel}

\usepackage{tikz}
\usetikzlibrary{cd}

\usepackage{xfrac}
\usepackage{todonotes}
\usepackage{wrapfig}

\usepackage{amsmath}
\usepackage{stmaryrd}
\usepackage{bm}

\usepackage{inqSem}

\newtheorem{PHONY}{PHONY}
\newtheorem{lemma}[PHONY]{Lemma}
\newtheorem{theorem}[PHONY]{Theorem}

\newtheorem{proposition}[PHONY]{Proposition}
\newtheorem{corollary}[PHONY]{Corollary}

\theoremstyle{definition}

\theoremstyle{definition}
\newtheorem{definition}[PHONY]{Definition}

\newcommand{\ol}[1]{\overline{#1}}
\newcommand{\tuple}[1]{\left\langle #1 \right\rangle}
\newcommand{\set}[1]{\left\{ #1 \right\}}

\newcommand{\N}{\mathbb{N}}

\newcommand{\formulas}{\mathcal{L}}
\newcommand{\sem}[1]{\left\llbracket #1 \right\rrbracket}
\newcommand{\class}[1]{\mathcal{#1}}

\newcommand{\Var}{\mathsf{Var}}
\newcommand{\Log}{\mathsf{Log}}

\newcommand{\IL}{\mathbf{IL}}
\newcommand{\HA}{\mathbf{HA}}
\newcommand{\ha}{\mathsf{HA}}

\newcommand{\Hvar}{\mathrm{H}}
\newcommand{\Svar}{\mathrm{S}}
\newcommand{\Pvar}{\mathrm{P}}

\newcommand{\latleq}{\preceq}

\newcommand\at{\mathrm{AP}}

\newcommand{\MP}{\mathrm{MP}}
\newcommand{\US}{\mathrm{US}}
\newcommand{\AP}{\mathrm{AP}}

\newcommand{\IPC}{\mathtt{IPC}}
\newcommand{\CPC}{\mathtt{CPC}}
\newcommand{\ND}{\mathtt{ND}}
\newcommand{\KP}{\mathtt{KP}}
\newcommand{\ML}{\mathtt{ML}}
\newcommand{\LC}{\mathtt{LC}}
\newcommand{\WEM}{\mathtt{WEM}}

\newcommand\vari{\mathcal{V}}

\newcommand\nvari{\mathcal{V}^\chi}

\newcommand\dvari{\mathcal{X}}

\tikzstyle{dot}=[inner sep=1pt, fill, black, circle, draw, minimum size = 5pt]
\tikzstyle{highlight}=[inner sep=1pt, draw=black, fill=white, circle, minimum size = 9pt]

\title{Lattices of Intermediate Theories via\\Ruitenburg's Theorem\thanks{The authors would like to thank Nick Bezhanishvili for comments and discussions on this work. The first author was supported by the European Research Council (ERC, grant agreement number 680220). The second author was supported by Research Funds of the University of Helsinki.}}

\author{
	Gianluca Grilletti
\and
    Davide Emilio Quadrellaro
}

\institute{
  Institute for Logic, Language and Computation (ILLC),
  Amsterdam, The Netherlands\\  
  Department of Mathematics and Statistics, University of Helsinki, Finland\\
  \email{grilletti.gianluca@gmail.com}\hspace{3em}
  \email{davide.quadrellaro@gmail.com}
  }

\authorrunning{Gianluca Grilletti, Davide Emilio Quadrellaro}

\titlerunning{Lattices of Intermediate Theories via Ruitenburg's Theorem}

\begin{document}

\maketitle

\begin{abstract}
	\noindent
	For every univariate formula $\chi$ we introduce a lattices of intermediate theories: the lattice of $\chi$-logics.
	The key idea to define $\chi$-logics is to interpret atomic propositions as fixpoints of the formula $\chi^2$, which can be characterised syntactically using Ruitenburg's theorem.
	We develop an algebraic duality between the lattice of $\chi$-logics and a special class of varieties of Heyting algebras.
	This approach allows us to build five distinct lattices---corresponding to the possible fixpoints of univariate formulas---among which the lattice of negative variants of intermediate logics.
	We describe these lattices in more detail.
\end{abstract}

\section{Introduction}

This paper introduces a family of lattices of intermediate theories, building on three results from the literature:
the duality between intermediate logics and varieties of Heyting algebras,
a novel algebraic semantics for inquisitive logic and negative variants,
and Ruitenburg's Theorem.

\smallskip
Intermediate logics (\cite{Chagrov:97,Gabbay:81}) are classes of formulas closed under uniform substitution and modus ponens, lying between intuitionistic logic ($\IPC$) and classical logic ($\CPC$).
This family of logics has been studied using several semantics, as for example Kripke semantics, Beth semantics, topological semantics and algebraic semantics (for an overview see \cite{BEZHANISHVILI2019403}).
Among these, the algebraic semantics based on Heyting algebras plays a special role:
every intermediate logic is sound and complete with respect to some class of Heyting algebras.%
\footnote{Kripke semantics is known to be incomplete for some intermediate logics. However, it is still an open problem whether this hold for Beth and topological semantics (\cite{Kuz75,BEZHANISHVILI2019403}).}

This connection between intermediate logics and Heyting algebras has been studied using tools from universal algebra.
As a consequence of Birkhoff's Theorem (\cite{Birkhoff1935,BurrisSankappanavar1981}), the lattice of varieties of Heyting algebras $\textbf{HA}$ is dually isomorphic to the lattice of intermediate logics $\textbf{IL}$.
This result allows to characterise properties of intermediate logics in terms of properties of the corresponding variety, and viceversa.

\smallskip
Inquisitive logic $\inqB$ (\cite{CiardelliRoelofsen:11jpl,Roelofsen:11,Ciardelli:16,Ciardelli:19book}) is an extension of classical logic that encompasses logical relation between \emph{questions} in addition to statements.
The logic was originally defined through the \emph{support semantics}, a generalisation of the standard truth-based semantics of $\CPC$.
Ciardelli et al. gave an axiomatisation of the logic, showing that it sits between $\IPC$ and $\CPC$, and highlighting connections with other intermediate logics such as Maksimova's logic $\ND$, Kreisel-Putnam logic $\KP$ and Medvedev's logic $\ML$ (\cite{Ciardelli:09thesis}).
However, $\inqB$ itself is not an intermediate logic, since it is not closed under uniform substitution.

An algebraic semantics for $\inqB$ has been defined in \cite{BezhanishviliGrillettiHolliday2019}, based on the corresponding algebraic semantics for intermediate logics.
The key idea is to restrict the interpretation of atomic propositions to range over regular elements of a Heyting algebra, that is, over fixpoints of the operator $\neg\neg$.
This restriction allows to have a sound and complete algebraic semantics, despite the failure of the uniform substitution principle.
As shown in \cite{Quadrellaro.2019}, this approach can be extended to the class of $\mathtt{DNA}$-logics, also known as \emph{negative variants of intermediate logics} (\cite{Miglioli:89,Iemhoff:2016aa}).
Moreover, this leads naturally to a duality between $\mathtt{DNA}$-logics and a special class of varieties, analogous to the one for intermediate logics.

\smallskip
Ruitenburg's theorem (\cite{Ruitenburg1984}) concerns sequences of formulas of the following kind:
\begin{equation*}
	\alpha^0  :=  p   \hspace{5em}
	\alpha^{n+1}  :=  \alpha[\sfrac{\alpha^n}{p}]
\end{equation*}
where $\alpha$ is a formula and $p$ is a fixed atomic proposition.
In particular, Ruitenburg's theorem states that this sequence is definitively periodic with period 2---modulo logical equivalence.
For example, if we take $\alpha := \neg p$ we can see that $\neg p \equiv \neg\neg\neg p$, showing that $\neg p$ is a fixpoint of the operator $\neg\neg$.
Ghilardi and Santocanale give an alternative proof of this result in \cite{santocanale:hal-01766636}, studying endomorphisms of finitely generated Heyting algebras.
This proof makes use of the duality introduced above and it highlights the relevance of the algebraic interpretation of Ruitenburg's Theorem.

\smallskip
In this paper we use Ruitenburg's theorem and its algebraic interpretation to define a lattice of intermediate theories in the same spirit as the negative variants.
Fixed a univariate formula $\chi$, we define an algebraic semantics by restricting valuations to range over fixpoints of the formula $\chi^2$---which can be characterised using Ruitenburg's Theorem.
This allows us to build the lattice of $\chi$-logics, intermediate theories characterised in terms of the fixpoint-axiom $\chi^2(p) \leftrightarrow p$.
We show that the algebraic semantics is sound and complete for these logics;
and we developed a duality theory for these logics analogous to the one for negative variants.
We also show that there are only six possible fixpoints for univariate formulas: $\top, p, \neg p, \neg\neg p, p\vee \neg p, \bot$.
This allows us to characterise and describe all the possible lattices of $\chi$-logics built using this approach.

In Section \ref{prel} we introduce some preliminary notions on intermediate logics and their algebraic semantics, the theory of algebraic duality for such logics and Ruitenburg's Theorem.
In Section \ref{chi} we define $\chi$-logics and give a brief overview of their main properties that can be derived in purely syntactic terms.
In Section \ref{alge}, fixed a formula $\chi$, we introduce a novel algebraic semantics for $\chi$-logics based on Ruitenburg's Theorem and we define a notion of variety of Heyting algebras suitable to study $\chi$-logics, namely $\chi$-varieties.
In Section \ref{dual} we develop an algebraic duality theory for these logics, showing that the lattice of $\chi$-logics is dually isomorphic to the lattice of $\chi$-varieties.
Finally, in Section \ref{discu} we show there are only 5 distinct lattices of $\chi$-logics for any univariate formula $\chi$, we describe their properties in more detail and we study the relations between them.
Conclusions and possible directions for future work are presented in Section \ref{conclusion}.

\section{Preliminaries}\label{prel}

In this Section we summarise the theory from the literature used throughout the paper.

\subsection{Algebraic semantics for intermediate logics}\label{subsection:algebraicSemantics}

Fix an infinite set $\AP$ of atomic propositions and consider the set of formulas $\formulas$ generated by the following grammar:
\begin{equation*}
	\phi \;::=\; p \;|\; \bot \;|\; \phi \land \phi \;|\; \phi \vee \phi \;|\; \phi \to \phi
\end{equation*}
where $p \in \AP$.
As usual, we will introduce the shorthand $\neg \phi := \phi \to \bot$ for \emph{negation}.
Henceforth we will leave the sets $\AP$ and $\formulas$ implicit, referring to \emph{atomic propositions from $\AP$} and to \emph{formulas from $\formulas$} simply as \emph{atomic propositions} and \emph{formulas} respectively.
To indicate a sequence of propositions $\tuple{p_1,\dots,p_n}$ we will often use the notation $\ol{p}$, and similarly for sequences of formulas $(\ol{\phi} = \tuple{\phi_1,\dots,\phi_n})$ and sequences of other objects.

Consider formulas $\phi, \psi$ and an atomic proposition $p$.
We will indicate with $\phi\left[ \sfrac{\psi}{p} \right]$ the formula obtained by substituting \emph{every} occurrence of $p$ in $\phi$ with the formula $\psi$.
More generally, given $\ol{\psi} = \tuple{\psi_1,\dots,\psi_n}$ a sequence of formulas and $\ol{p} = \tuple{p_1,\dots,p_n}$ a sequence of \emph{distinct} atomic propositions, we will indicate with $\phi\left[ \sfrac{\ol{\psi}}{\ol{p}} \right]$ the formula obtained by substituting \emph{simultaneously} each $p_i$ with $\psi_i$.
With abuse of notation, when we take $\chi$ a univariate formula---that is, a formula with only one free variable---we will indicate the sequence $\tuple{\chi(p_1),\dots,\chi(p_n)}$ with the notation $\chi(\ol{p})$;
for example, the notations $\phi[\sfrac{\neg\ol{p}}{\ol{p}}]$ and $\phi[\sfrac{\tuple{\neg p_1,\dots,\neg p_n}}{\tuple{p_1,\dots,p_n}}]$ indicate the same formula.



\medskip
We will refer to the \emph{intuitionistic propositional calculus} (see for example \cite{Chagrov:97,Gabbay:81}) as $\IPC$.
With slight abuse of notation, we will write $\IPC$ also to refer to the set of validities of this calculus.

An \emph{intermediate logic} (\cite{Chagrov:97,Gabbay:81}) is a set of formulas $L$ with the following properties:
\begin{enumerate}
	\item $\IPC \subseteq L \subseteq \CPC$;
	\item $L$ is closed under \emph{modus ponens}: If $\phi \in L$ and $\phi \to \psi \in L$, then $\psi \in L$;
	\item $L$ is closed under \emph{uniform substitution}: If $\phi \in L$ and $\rho$ is a substitution, then $\phi[\rho] \in L$.
\end{enumerate}

\noindent
Given $\Gamma$ a set of formulas, we indicate with $\MP(\Gamma)$ the smallest set of formulas extending $\Gamma$ and closed under modus ponens;
and with $\US(\Gamma)$ the smallest set of formulas extending $\Gamma$ and closed under uniform substitution.
It is immediate to prove such sets always exist and that if $\Gamma \subseteq \CPC$, then $\MP(\US(\Gamma))$ is the smallest intermediate logic extending $\Gamma$:
we will call $\MP(\US(\Gamma))$ the intermediate logic \emph{generated} by $\Gamma$.

Intermediate logics form a structure of \emph{bounded lattice} under the set-theoretic inclusion.
In particular, the meet and join operations are $L_1 \land L_2 := L_1 \cap L_2$ and $L_1 \vee L_2 := \MP(L_1 \cup L_2)$, and $\IPC$ and $\CPC$ are the minimum and maximum respectively.
We will refer to this lattice with the notation $\IL$.

\medskip
In the literature, several semantics have been proposed to study these logics:
Kripke semantics, Beth semantics and topological semantics are some famous examples (see \cite{BEZHANISHVILI2019403} for an overview of some well-known semantics).
In this paper we will focus on the so-called \emph{algebraic semantics}:
given an \emph{Heyting algebra} $H$ and a function $V: \AP \to H$---which we will refer to as a \emph{valuation}---we define recursively by the following clauses the interpretation $\llbracket \phi \rrbracket^{H}_V$ of a formula $\phi$ in $H$ under $V$.
\begin{equation*}
\begin{array}{r@{\hspace{.3em}}c@{\hspace{.3em}}l  @{\hspace{1.5em}}  r@{\hspace{.3em}}c@{\hspace{.3em}}l  @{\hspace{1.5em}}  r@{\hspace{.3em}}c@{\hspace{.3em}}l}
	\llbracket p \rrbracket^{H}_{V} &= &V(p)
	&\llbracket \top \rrbracket^{H}_{V} &= &1_H
	&\llbracket \bot \rrbracket^{H}_{V} &= &0_H \\[.5em]
	\llbracket \phi \land \psi \rrbracket^{H}_{V} &= &\llbracket \phi \rrbracket^{H}_{V} \land_H \llbracket \psi \rrbracket^{H}_{V}
	&\llbracket \phi \vee \psi \rrbracket^{H}_{V} &= &\llbracket \phi \rrbracket^{H}_{V} \vee_H \llbracket \psi \rrbracket^{H}_{V}
	&\llbracket \phi \to \psi \rrbracket^{H}_{V} &= &\llbracket \phi \rrbracket^{H}_{V} \to_H \llbracket \psi \rrbracket^{H}_{V}
\end{array}
\end{equation*}
where $1_H,0_H,\land_H,\vee_H,\to_H$ indicate the constants and operations of the algebra $H$.
We say that $\phi$ is \textit{true in $H$ under $V$} and we write $(H,V)\vDash \phi$ if $\sem{\phi}^H_V = 1$.
We say that $\phi$ is \textit{valid in $H$} and we write $H\vDash \phi$ if it is true in $H$ under any valuation $V$.

As a shorthand, we will indicate with $[p_1\mapsto a_1,\dots,p_n\mapsto a_n]$ a generic valuation $V$ such that $V(p_i) = a_i$---without specifying its value on the atomic formulas different from $p_1,\dots, p_n$.

Given a function $f: H^n \to H$ we will call it a \emph{polynomial} if it is obtained by composing the functions $1_H, 0_H, \land_H, \vee_H$ and $\to_H$ (where we identify the constants $1_H$ and $0_h$ with the corresponding $0$-ary functions).
Given a formula $\phi(p_1,\dots,p_n)$, we can associate to it the polynomial $\bm{\phi}$ (indicated with the bold font) defined as:
\begin{equation*}
\begin{array}{rlcl}
	\bm{\phi}: &H^n &\to &H \\
		&\ol{a} &\mapsto &\llbracket \phi \rrbracket^{H}_{[p_1\mapsto a_1,\dots,p_n\mapsto a_n]}
\end{array}
\end{equation*}
Moreover, it is immediate to show that for every polynomial $f$, there exists a (non-unique) formula $\phi$ such that $f = \bm{\phi}$.

\subsection{Algebraic duality}

In the study of intermediate logics, a special role is played by the collection of those algebras \emph{defined} by a certain intermediate logic.
Given $L \in \IL$, define the \emph{variety generated by $L$} as the set
\begin{equation*}
	\Var(L) \;:=\; \{\;  H \in \ha \;|\; \forall \phi \in L.\, H\vDash \phi \;\}
\end{equation*}
where $\ha$ indicates the class of all Heyting algebras.
We will call a class $\vari \subseteq \ha$ a \emph{variety} if $\vari$ is closed under the operations $\Hvar, \Svar, \Pvar$ defined over subclasses of $\ha$ as follows:
\begin{align*}
	\Hvar(\class{C}) &:= \{\;  H \in \ha \;|\;  \exists A \in \class{C}.\, A \twoheadrightarrow H  \;\}  &&\text{(homomorphic images)}\\
	\Svar(\class{C}) &:= \{\;  H \in \ha \;|\;  \exists A \in \class{C}.\, H \hookrightarrow A  \;\}  &&\text{(subalgebras)}\\
	\Pvar(\class{C}) &:= \left\{\;  \prod_{i\in I} A_i \in \ha \;\middle|\;  \forall i\in I.\, A_i \in \class{C}  \;\right\}  &&\text{(products)}
\end{align*}

\noindent
It is easy to prove that $\Var(L)$ is indeed a variety;
moreover the following well-known results give us an alternative characterisation of varieties:

\begin{theorem}[Tarski's theorem; \cite{Tarski1946}, Theorem 9.5 in \cite{BurrisSankappanavar1981}]
	Given $\class{C} \subseteq \ha$ a class of algebras, $\bm{\vari}(\class{C}) := \Hvar\Svar\Pvar(\class{C})$ is the smallest variety containing $\class{C}$.
\end{theorem}

\begin{theorem}[Birkhoff's theorem; \cite{Birkhoff1935}, Theorem 11.9 in \cite{BurrisSankappanavar1981}]
	A class of algebras $\class{C} \subseteq \ha$ is a variety iff it is equationally definable, that is, there exists a set of formulas $F \subseteq \mathcal{L}$ such that
	\begin{equation*}
		\class{C} = \{\;  H \in \ha  \;|\;  \forall \phi \in F.\, H\vDash \phi  \;\}.
	\end{equation*}
\end{theorem}

\noindent
Given these results, it is easy to show that varieties form a \emph{bounded distributive lattice} under the inclusion order.
In particular, the meet and join operations are $\vari_1 \land \vari_2 := \vari_1 \cap \vari_2$ and $\vari_1 \vee \vari_2 := \bm{\vari}(\vari_1 \cup \vari_2)$.
We will refer to this lattice with the notation $\HA$.\footnote{Notice the difference between $\ha$ (the class of all Heyting algebras) and $\HA$ (the lattice of \emph{varieties} of Heyting algebras). In particular $\ha\in \HA$.}

Moreover, given a variety $\vari$ we can define a set of formulas that characterises it:
\begin{equation*}
	\Log(\vari) := \{\; \phi \in \mathcal{L} \;|\; \forall H\in \vari.\, H\vDash \phi \;\}
\end{equation*}
It is easy to prove that $\Log(\vari)$ is an intermediate logic;
and that---using Birkhoff's theorem---$L = \Log(\Var(L))$ and $\vari = \Var(\Log(\vari))$ for every intermediate logic $L$ and variety $\vari$.
Moreover, since $\Var$ and $\Log$ are antitone maps, these maps are dual isomorphisms between the lattice $\IL$ and the lattice $\HA$:

\begin{theorem}[Duality; Theorem 7.54 in \cite{Chagrov:97}]
	The lattice of intermediate logics is dually isomorphic to the lattice of varieties of Heyting algebras, i.e. $\IL \cong^{op}\HA$.
\end{theorem}

\subsection{Ruitenburg's Theorem}


For the remainder of this Section, we will indicate with $p$ a fixed atomic proposition.
Let $\phi(p,\ol{q})$ be a formula, where $p,\ol{q}$ contain all the atomic propositions appearing in $\phi$.
A folklore result says that the formulas
\begin{equation*}
	\phi(p,\ol{q}) \hspace{5em} \phi^{3}(p,\ol{q}) := \phi(\,\phi(\,\phi(p,\ol{q}),\,\ol{q}),\,\ol{q})
\end{equation*}
are equivalent in classical logic.
Ruitenburg extended this result to intuitionistic logic (\cite{Ruitenburg1984}).

\begin{definition}
	Given $\phi(p,\ol{q})$ a formula, define the formulas $\{ \phi^n(p,\ol{q}) \}_{n\in \N}$ recursively as follows:
	\begin{equation*}
		\phi^0(p,\ol{q}) := p  \hspace{5em}  \phi^n(p,\ol{q}) := \phi(\, \phi^{n-1}(p,\ol{q}),\, \ol{q}\, )
	\end{equation*}
	That is, $\phi^n$ is obtained by substituting $\phi^{n-1}$ for $p$ in $\phi$.
\end{definition}

\begin{theorem}[Ruitenburg's theorem; \cite{Ruitenburg1984}]\label{theorem:ruitenburg}
	For every formula $\phi(p,\ol{q})$, the sequence $\phi^0,\phi^1,\phi^2,\dots$ is---modulo logical equivalence---definitely periodic with period $2$.
	That is, there exists a natural number $n$ such that:
	\begin{equation}\label{eq:index}
		\phi^n \leftrightarrow \phi^{n+2}  \in  \IPC
	\end{equation}
\end{theorem}

\noindent
We will call the smallest $n$ for which Equation \ref{eq:index} holds the \emph{Ruitenburg index} (or simply the \emph{index}) of $\phi$.
Moreover, we will call $\phi^n$ the \emph{fixpoint} of the formula $\phi$.

We can also see Ruitenburg's result as an algebraic fixpoint theorem.
Let $A$ be a Heyting algebra, $\ol{a}$ a sequence of elements in $A$ and $f(x,\ol{y})$ a polynomial.
Then a consequence of Ruitenburg's theorem is that the operator $f^2(x,\ol{a}) = f(f(x,\ol{a}),\ol{a})$ admits a fixpoint.
And indeed, this is an equivalent formulation of Theorem \ref{theorem:ruitenburg}, as can be easily shown by applying it to the Lindebaum-Tarski algebra of $\IPC$.

As proven by Ruitenburg (Example 2.5 in \cite{Ruitenburg1984}), there is no uniform bound for all the formulas $\phi$, but each formula admits an index.
However, for some classes of formulas we can find a uniform bound:

\begin{lemma}[Proposition 2.3 in \cite{Ruitenburg1984}]\label{lemma:finiteFixpoints}
	If $\chi(p)$ is a univariate formula, then
		$\chi^2 \leftrightarrow \chi^4 \in \IPC.$
	Moreover, the fixpoint of $\chi(p)$ is equivalent to one of the following formulas: $\bot$, $p$, $\neg p$, $\neg\neg p$, $p \vee \neg p$, $\top$.
\end{lemma}

\noindent
We give an elementary proof of this result in Appendix \ref{appendixA}, different from the original one given by Ruitenburg in \cite{Ruitenburg1984}.

\section{$\chi$-logics}\label{chi}

In analogy with the case of negation, given a formula $\chi$ we are interested in logics arising by interpreting atoms as fixpoints for the operator $\chi^2$:
we will call these logics \emph{$\chi$-logics}.

In this paper, we start the study of this family of logics by considering only univariate formulas for two reasons:
the presence of additional atoms requires a generalisation of the duality results presented in Section \ref{prel} to Heyting algebras with constants;
and restricting our attention to univariate formulas allows us to give a more in-depth description of all the lattices of $\chi$-logics generated through this procedure---which are finitely many, as shown at the end of this section.

\begin{definition}[$\chi$-logic]\label{definition:chiLogic}
	Let $\chi(p)$ be a univariate formula and $\Gamma$ a set of formulas.
	We define the \emph{$\chi$-logic generated by $\Gamma$} as the smallest set of formulas $\Gamma^{\chi}$ with the following properties:
	\begin{enumerate}
		\item $\IPC \subseteq \Gamma^{\chi}$; \label{condition:ipc}
		\item If $\phi \in \Gamma$ and $\sigma$ is a substitution, then $\phi[\sigma] \in \Gamma^{\chi}$; \label{condition:gammaUS}
		\item $\chi^{2}(p) \leftrightarrow p \in \Gamma^{\chi}$ for every atomic proposition $p$; \label{condition:fixpoint}
		\item $\Gamma^{\chi}$ is closed under \emph{modus ponens}: if $\phi \in \Gamma^{\chi}$ and $\phi \to \psi \in \Gamma^{\chi}$, then $\psi \in \Gamma^{\chi}$. \label{condition:modusPonens}
	\end{enumerate}
\end{definition}

\noindent
Condition \ref{condition:fixpoint} requires atoms to behave as fixpoints of the operator $\chi^2$.
This condition together with uniform substitution would impose that \emph{all} formulas behave like fixpoints, which is a requirement too strong for our purposes (we will see later that $\chi$-logics are generally \emph{not} closed under uniform substitution).
That is why we require the uniform substitution principle only for formulas in $\Gamma$, that is, Condition \ref{condition:gammaUS}.

Notice that we can interpret $\Gamma^{\chi}$ as the set of valid formulas of a Hilbert-style deductive system where Conditions \ref{condition:ipc}, \ref{condition:gammaUS} and \ref{condition:fixpoint} define the axioms---to be more precise, the axiom schemata---while Conditions \ref{condition:modusPonens} specifies modus ponens as the only rule of the system.
This suggests the following characterisation of $\chi$-logics.

\begin{lemma}\label{lemma:gammaOrL}
	Let $L$ be the intermediate logic generated by $\Gamma$.
	Then $\Gamma^{\chi} = L^{\chi}$.
\end{lemma}
\begin{proof}
	The left-to-right containment is immediate, since the operator $(-)^{\chi}$ is monotone.
	As for the other containment, notice that Conditions $\ref{condition:ipc}$, \ref{condition:gammaUS} and \ref{condition:modusPonens} impose that $L \subseteq \Gamma^{\chi}$, from which the result follows.
\end{proof}

\noindent
So we can think of $\chi$-logics as always generated by a corresponding intermediate logic instead of a generic set of formulas.
To stress this point, given an intermediate logic $L$ we will call $L^{\chi}$ the \emph{$\chi$-variant of $L$}.
Notice that a direct consequence of Lemma \ref{lemma:gammaOrL} is that any set satisfying Conditions \ref{condition:gammaUS}, \ref{condition:fixpoint} and \ref{condition:modusPonens} is the $\chi$-variant of some intermediate logic $L$.

Restricting our attention to intermediate logics allows us to give an alternative characterisation of $\chi$-logics.

\begin{lemma}\label{lemma:characterizationChiVariant}
	Let $\chi(p)$ be a univariate formula and $n$ be its index.
	Given $L$ an intermediate logic, we have
	\begin{equation*}
		L^{\chi} \;=\; \set{\;  \phi(\ol{p})  \;\middle|\;  \phi[\sfrac{\chi^{n}(\ol{p})}{\ol{p}}] \in L  \;}.
	\end{equation*}
\end{lemma}

\noindent

\begin{proof}
	Call the set on the right-hand side $M$.

	Firstly, we will show that $M$ satisfies the conditions in Definition \ref{definition:chiLogic}.
	Since $L$ contains $\IPC$ and is closed under modus ponens and uniform substitution, we easily obtain Conditions \ref{condition:ipc}, \ref{condition:gammaUS} and \ref{condition:modusPonens}.
	As for Condition \ref{condition:fixpoint}, since $n$ is the index of $\chi$, we have $\chi^{n+2}(p) \leftrightarrow \chi^{n}(p) \in L$, from which it follows $\chi^2(p) \leftrightarrow p \in M$ for every atomic proposition $p$.

	Secondly, we need to show that $M$ is the smallest set satisfying these conditions.
	To do so, we will use the following fact:
	given formulas $\ol{\alpha} = \tuple{\alpha_1, \dots, \alpha_l}, \ol{\beta} = \tuple{ \beta_1, \dots, \beta_l}, \gamma$ formulas and distinct atomic propositions $\ol{q} = \tuple{q_1,\dots,q_l}$, we have
	\begin{equation*}
		\bigwedge_{i\leq l}(\; \alpha_i \leftrightarrow \beta_i \;) \rightarrow (\; \gamma[\sfrac{\ol{\alpha}}{\ol{q}}] \leftrightarrow \gamma[\sfrac{\ol{\beta}}{\ol{q}}] \;) \in \IPC \subseteq X
	\end{equation*}

	\noindent
	Consider now a set $X$ satisfying the conditions of Definition \ref{definition:chiLogic}.

	Since $\chi^2(q) \leftrightarrow q \in X$ for every $q$ and $X$ is closed under uniform substitution, it follows that also $\chi^{4}(q) \leftrightarrow \chi^{2}(q)\in X$.
	Moreover, since
	\begin{equation*}
		(\alpha \leftrightarrow \beta) \rightarrow (\; (\beta \leftrightarrow \gamma) \rightarrow (\alpha \leftrightarrow \gamma) \;) \in \IPC \subseteq X
	\end{equation*}
	by closure under modus ponens we obtain that $\chi^4(q) \leftrightarrow q \in X$.
	Iterating this reasoning, we obtain that $\chi^n(q) \leftrightarrow q \in X$ for every $q$; or $\chi^{n+1}(q) \leftrightarrow q$ for every $q$---depending on the parity of $n$.
	Assume the former is the case;
	the treatment of the other case is analogous.

	Consider now a formula $\phi(\ol{p})$ with $\ol{p} = \tuple{p_1,\dots,p_l}$. Combining the previous facts we get:
	\begin{equation*}
	\begin{array}{ccl}
		&\chi^{n}(p_i) \leftrightarrow p_i &\in X \;\text{for every $i\leq l$}\\
		\text{and} &\bigwedge_{i\leq l}(\; \chi^{n}(p_i) \leftrightarrow p_i \;) \rightarrow (\; \phi[\sfrac{\chi^{n}(\ol{p})}{\ol{p}}] \leftrightarrow \phi(\ol{p}) \;) &\in X   \\[.5em]
		\text{implies}  &\phi[\sfrac{\chi^{n}(\ol{p})}{\ol{p}}] \leftrightarrow \phi(\ol{p}) &\in X
	\end{array}
	\end{equation*}

	\noindent
	Suppose now $\phi[\sfrac{\chi^n(\ol{p})}{\ol{p}}] \in L$.
	Since
	\begin{equation*}
		\alpha \rightarrow (\;  (\alpha \leftrightarrow \beta) \rightarrow \beta  \;) \in \IPC \subseteq X
	\end{equation*}
	and $X$ is closed under modus ponens, it follows that $\phi(\ol{p}) \in X$; and consequently $M \subseteq X$.

	So $M$ is the smallest set satisfying the Conditions in Definition \ref{definition:chiLogic}, thus proving $L^{\chi} = M$ as wanted.
\end{proof}

\noindent
It is interesting to notice that instances of $\chi$-logics
 have already been studied in the literature:
an example is \emph{inquisitive logic} $\inqB$.
In fact, reinterpreting Theorem 3.4.9 in \cite{Ciardelli:09thesis}, we obtain the following:

\begin{theorem}[Ciardelli; Theorem 3.4.9 in \cite{Ciardelli:09thesis}]\label{inq}
	\begin{equation*}
		\KP^{\neg p} = \ND^{\neg p} = \ML^{\neg p} = \inqB
	\end{equation*}
\end{theorem}

\noindent
One can easily show that the set $L^\chi$ in general is not closed under \emph{uniform substitution}.
Nonetheless, the set $L^{\chi}$ is closed under a weaker notions of substitution, that is, \emph{atomic substitution}:
give $\sigma \in \AP \to \AP$ a permutation of the atomic propositions, if $\phi(\ol{p}) \in L^{\chi}$ then $\phi(\sigma(\ol{p}))$.
That is to say, even though atomic propositions play a special role---fix-points of $\chi^2$---they are still considered as generic entities, in that they are indistinguishable from one another.
Moreover, as noted by Iemhoff and Yang in \cite{Iemhoff:2016aa}, the logics $L^{\neg p}$ are closed under a more general substitution principle, that is, \emph{classical substitutions}:
A classical substitution maps every atomic proposition with a $\vee$-free formula.
In general, it is expected that for a fixed $\chi$ the logics $L^{\chi}$ are closed under more general substitution principles.




We can show that $\chi$-logics form a bounded distributive lattice under the set-theoretic containment, as it was the case for intermediate logics.
In particular, the meet operation is given by set-theoretic intersection and the join operation by the closure under modus ponens of the union---in complete analogy with the case of intermediate logics.

\begin{lemma}
	Given $\chi$ a univariate formula and $L, M$ two intermediate logics we have:
	\begin{equation*}
		L^{\chi} \land M^{\chi} = L^{\chi} \cap M^{\chi} = (L \land M)^{\chi}  \qquad\qquad
		L^{\chi} \vee M^{\chi} = \MP( L^{\chi} \cup M^{\chi}) = (L \vee M)^{\chi}
	\end{equation*}
\end{lemma}
\begin{proof}
	We consider only the second set of identities, as the proof can be easily adapted for the first set.
	Firstly, notice that $L^{\chi}, M^{\chi} \subseteq (L\vee M)^{\chi}$.
	Moreover, since $L\subseteq L^{\chi}$ and $M\subseteq M^{\chi}$, for every $\chi$-logic $\Lambda$ such that $L^{\chi},M^{\chi} \subseteq \Lambda$ it holds $L\cup M \subseteq \Lambda$;
	and since $\chi$-logics are closed under modus ponens it holds $L \vee M = \MP(L \cup M) \subseteq \Lambda$.
	So in particular $(L\vee M)^{\chi} \subseteq \Lambda$.
	This implies that $(L\vee M)^{\chi}$ is the least upper bound of $L^{\chi}$ and $M^{\chi}$, that is, $L^{\chi} \vee M^{\chi} = (L\vee M)^{\chi}$.

	Secondly, notice that $\MP(L^{\chi}\cup M^{\chi})$ is the $\chi$-logic generated by the set of formulas $L^{\chi}\cup M^{\chi}$.
	So in particular, since $L^{\chi},M^{\chi} \subseteq \MP(L^{\chi}\cup M^{\chi})$, we also have $L^{\chi} \vee M^{\chi} \subseteq \MP(L^{\chi}\cup M^{\chi})$.
	Moreover, since $(L\vee M)^{\chi}$ is closed under modus ponens and $L^{\chi} \cup M^{\chi} \subseteq (L\vee M)^{\chi}$, it follows $\MP(L^{\chi} \cup M^{\chi}) \subseteq (L\vee M)^{\chi} = L^{\chi} \vee M^{\chi}$.
	From this we conclude that $L^{\chi} \vee M^{\chi} = \MP( L^{\chi} \cup M^{\chi})$, as wanted.


\end{proof}

We will indicate with $\IL^{\chi}$ the lattice of $\chi$-logics.
Notice that the previous proof shows also that the mapping $L \mapsto L^{\chi}$ is a lattice morphism.

In the next sections we will study the structure of this lattice employing some tools from algebraic semantics.
But before moving to that, we will tackle one last question in this Section:
how many lattices are we dealing with?

As noted in Lemma \ref{lemma:finiteFixpoints}, there are only a finite amount of fixpoints, from which the following result follows readily.

\begin{theorem}\label{theorem:finitelyManyFixpoints}
	There are only 6 Ruitenburg-fixpoints of univariate intuitionistic formulas: $\bot$, $p$, $\neg p$, $\neg\neg p$, $p \vee \neg p$ and $\top$.
\end{theorem}

\noindent
Notice also that $\IL^{\neg p} = \IL^{\neg\neg p}$ are the same lattice: this follows from Lemma \ref{lemma:characterizationChiVariant}, since
\begin{equation*}
	L^{\neg p}
	\;=\;  \{\;  \phi(\ol{q})  \;|\;  \phi[\sfrac{\neg\ol{q}}{\ol{q}}]  \;\}
	\;=\;  \{\;  \phi(\ol{q})  \;|\;  \phi[\sfrac{\neg\neg\ol{q}}{\ol{q}}]  \;\}
	\;=\;  L^{\neg\neg p}
\end{equation*}
for every intermediate logic $L$.
So in total we are working with only $5$ lattices, associated to the $6$ fix-points above.
In Section \ref{discu} we will see that these are indeed distinct lattices.
\begin{equation*}
	\IL^{\bot} \hspace{4em} \IL^{p} = \IL \hspace{4em} \IL^{\neg p} = \IL^{\neg\neg p} \hspace{4em} \IL^{p \vee \neg p} \hspace{4em} \IL^{\top}
\end{equation*}

\section{Algebraic semantics}\label{alge}

In this section we shall provide a semantic interpretation of $\chi$-logics, by generalizing the algebraic semantics for inquisitive logic presented in \cite{BezhanishviliGrillettiHolliday2019} and further developed in \cite{Quadrellaro.2019}.
The key to generalise the algebraic semantics to this context lies in an algebraic interpretation of Ruitenburg's Theorem.
In this section we will fix a univariate formula $\chi$ with index $n$.

\begin{wrapfigure}{r}{0.45\textwidth}
	\centering
	\begin{tikzpicture}
	\node at (-2.1,1.7) {$H$};
	\draw (0,0) ellipse (3cm and 2cm);
	
	\node at (-1.2,1.4) {$\chi[H]$};
	\draw (0,0) ellipse (2.2cm and 1.3cm);
	
	\node at (-1.6,0) {$\cdots$};
	
	\node at (0,.8) {$\chi^n[H]$};
	\draw[thick] (0,0) ellipse (1cm and .5cm);
	
	\node[circle,fill=black,inner sep=0cm,minimum size=.15cm] (a) at (-.6,0) {};
	\node[circle,fill=black,inner sep=0cm,minimum size=.15cm] (b) at (  0,0) {};
	
	\node[circle,fill=black,inner sep=0cm,minimum size=.15cm] (c) at (.4,-.2) {};
	
	\draw[-latex] (a) to[bend left] (b);
	\draw[-latex] (b) to[bend left] (a);
	
	\draw[-latex] (c) to[loop above] (c);
	
	\end{tikzpicture}
\end{wrapfigure}
	

As noted in Section \ref{prel}, given a Heyting algebra $H$ we can define a polynomial corresponding to $\chi$:
\begin{align*}
	\bm{\chi}: H\;\rightarrow\; &H \\
	a \;\mapsto\; &\sem{\chi(p)}_{[p\mapsto a]}
\end{align*}

\noindent
Ruitenburg's Theorem tells us that the sequence $H, \bm{\chi}[H], \bm{\chi}^{2}[H] := \bm{\chi}[\bm{\chi}[H]],\dots$ is definitely constant;
and that $\bm{\chi}$ restricted to the set $\bm{\chi}^n[H]$ is an \emph{involution}.
Henceforth we will call the set $H^{\chi} := \bm{\chi}^n[H]$ the \emph{$\chi$-core} (or simply \emph{core} when $\chi$ is clear from the context) of $H$.
Notice that the $\chi$-core consists exactly of the fixpoints of $\bm{\chi}^2$:

\begin{lemma}\label{lemma:fixpointsAlgebras}
	$H^{\chi}$ is the set of fixpoints of $\bm{\chi}^2$.
\end{lemma}
\begin{proof}
	By Theorem \ref{theorem:ruitenburg}, we have that $\chi^n \equiv \chi^{n+2}$.
	Consider now an element in $a \in H^{\chi}$, that is, an element of the form $a = b^n$ for some $b\in H$.
	It follows that
	\begin{equation*}
		\bm{\chi}^2(a) = \bm{\chi}^2(\bm{\chi}^n(b)) = \bm{\chi}^{n+2}(b) = \bm{\chi}^n(b) = a
	\end{equation*}
	showing that $a$ is a fixpoint of $\bm{\chi}^2$.
	Conversely, let $a$ be a fixpoint for $\bm{\chi}^2$.
	Then it follows that
	\begin{equation*}
		a \quad=\quad \bm{\chi}^2(a) \quad=\quad  \bm{\chi}^2(\bm{\chi}^2(a)) = \bm{\chi}^4(a)   \quad=\quad \dots \quad=\quad \bm{\chi}^{2n}(a) \in H^{\chi}.
	\end{equation*}
\end{proof}

\noindent
For instance, when $\chi(p)=\neg p$ the core $H^{\neg p}$ of $H$ consists of the \emph{regular elements} of the algebra $H$, that is, fixpoints of the operator $\neg\neg$.

To obtain an adequate semantics for a $\chi$-logic, it is sufficient to restrict the valuations of atomic propositions to the core $H^\chi$. Let $\at$ be an arbitrary set of atomic propositions. We say that a valuation $\sigma: \at \to H$ is a \emph{$\chi$-valuation} if $\sigma[\at]\subseteq H^\chi$. A $\chi$-valuation over $H$ thus sends every atomic proposition to some element of the $\chi$-core of $H$. Algebraic models of $\chi$-logics are then defined as follows.

\begin{definition}[$\chi$-Model]
	A \textit{$\chi$-model} is a pair $M=(H,\sigma)$ such that $H$ is a Heyting algebra and $\sigma$ is a $\chi$-valuation.
\end{definition}

\noindent
The \textit{interpretation} of a formula $\phi\in\mathcal{L}$ in $M=(H,\sigma)$, in symbols $\sem{\phi}^{H}_{\sigma}$, can then be easily defined recursively, as in the standard algebraic semantics over Heyting algebras:
\begin{equation*}
\begin{array}{r@{\hspace{.3em}}c@{\hspace{.3em}}l  @{\hspace{1.5em}}  r@{\hspace{.3em}}c@{\hspace{.3em}}l  @{\hspace{1.5em}}  r@{\hspace{.3em}}c@{\hspace{.3em}}l}
	\llbracket p \rrbracket^{H}_{\sigma} &= &\sigma(p)
	&\llbracket \top \rrbracket^{H}_{\sigma} &= &1_H
	&\llbracket \bot \rrbracket^{H}_{\sigma} &= &0_H \\[.5em]
	\llbracket \phi \land \psi \rrbracket^{H}_{\sigma} &= &\llbracket \phi \rrbracket^{H}_{\sigma} \land_H \llbracket \psi \rrbracket^{H}_{\sigma}
	&\llbracket \phi \vee \psi \rrbracket^{H}_{\sigma} &= &\llbracket \phi \rrbracket^{H}_{\sigma} \vee_H \llbracket \psi \rrbracket^{H}_{\sigma}
	&\llbracket \phi \to \psi \rrbracket^{H}_{\sigma} &= &\llbracket \phi \rrbracket^{H}_{\sigma} \to_H \llbracket \psi \rrbracket^{H}_{\sigma}
\end{array}
\end{equation*}
where $1_H,0_H,\land_H,\vee_H,\to_H$ indicate the constants and operations of the algebra $H$.

The key point is that for atomic propositions $\llbracket p \rrbracket^{H}_{\sigma}= \sigma(p) \in H^{\chi}$, which means that the interpretation of every atomic proposition is a fixpoint for $\chi^2$.
We say that $\phi$ is \textit{true in $H$ under $\sigma$} and we write $(H,\sigma)\vDash^\chi \phi$ if   $\sem{\phi}^H_{\sigma} = 1$.
We say that $\phi$ is \textit{valid in $H$} and we write   $H\vDash^\chi \phi$ if it is true in $H$ under any   $\chi$-valuation $\sigma$.

The algebraic semantics we have introduced differs from the standard semantics of intermediate logics in the fact that we consider only a particular class of valuations for the atomic propositions, namely $\chi$-valuations.
The relation between validity at a Heyting algebra and $\chi$-validity is further clarified by the following results.
We define, for every valuation $V:\at\rightarrow H$, its \emph{$\chi$-variant} $V^\chi$ as the $\chi$-valuation $V^\chi:\at\rightarrow H^\chi$ such that $V^\chi(p)=\chi^{n}(V(p))$.
With a simple induction we can show the following connection between $V$ and $V^{\chi}$:
$\llbracket \phi \rrbracket^H_{V^\chi} = \left\llbracket \phi\left[ \sfrac{\chi^n(\ol{p})}{\ol{p}} \right] \right\rrbracket^H_V$.




Notice that, since $H^{\chi}$ is the image of $\chi^n$, every $\chi$-valuation is the $\chi$-variant of some valuation.
In fact, given $\sigma$ a $\chi$-valuation, $\sigma(p) \in H^{\chi} = \chi^{n}[H]$.
So, for any valuation $V$ such that $V(p) \in (\chi^{n})^{-1}(\sigma(p))$, we have $V^\chi(p) = \chi^n(V(p)) = \sigma(p)$.
This allows us to prove the following Lemma connecting validity in the standard algebraic sense and in the context of $\chi$-logics:

\begin{proposition}\label{lem8.5}
	For any Heyting algebra $H$, $H\vDash^\chi \phi$ if and only if $H\vDash \phi\left[\sfrac{\chi^n(\ol{p})}{\ol{p}} \right] $.
\end{proposition}
\begin{proof}
	Assume $H\nvDash \phi\left[ \sfrac{\chi^n(\ol{p})}{\ol{p}} \right]$, that is, there exists a valuation $V$ such that $\left\llbracket \phi\left[ \sfrac{\chi^n(\ol{p})}{\ol{p}} \right] \right\rrbracket^H_V \neq 1_H$.
	Considering the $\chi$-valuation $V^{\chi}$, we then have $\llbracket\phi\rrbracket^H_{V^{\chi}} = \left\llbracket \phi\left[ \sfrac{\chi^n(\ol{p})}{\ol{p}} \right] \right\rrbracket^H_V \neq 1_H$, and thus $H \nvDash^{\chi} \phi$.

	Conversely, assume $H \nvDash^{\chi} \phi$, that is, there exists a $\chi$-valuation $\sigma$ such that $\llbracket\phi\rrbracket^H_{\sigma} \neq 1_H$.
	As noted above, for some valuation $V$ we have $\sigma = V^\chi$;
	from this we obtain $\left\llbracket \phi\left[ \sfrac{\chi^n(\ol{p})}{\ol{p}} \right] \right\rrbracket^{H}_{V} = \llbracket\phi\rrbracket^H_{V^\chi} \neq 1_H$, and thus $H \nvDash \phi\left[ \sfrac{\chi^n(\ol{p})}{\ol{p}} \right]$.
\end{proof}

\noindent
Combining Lemma \ref{lemma:characterizationChiVariant} with the previous Proposition, we obtain the following Corollary.

\begin{corollary}\label{lem9}
	Let $H$ be a Heyting algebra and $L$ an intermediate logic, if $H\vDash L$ then $H\vDash^\chi L^\chi$.
\end{corollary}

\noindent
The converse of the previous Corollary does not hold in general, as a formula might be true in a Heyting algebra under all $\chi$-valuations but not under all valuations.
The next proposition is a weaker version of this converse:
Let $\langle H^\chi \rangle$ be the subalgebra of $H$ generated by the $\chi$-core $H^\chi$; we say that $H$ is \textit{core generated} if $H=\langle H^\chi \rangle$.

\begin{lemma}\label{lem9.5}
	Let  $H$ be a Heyting algebra,  $H\vDash^\chi \phi $ if and only if $\langle H^\chi \rangle\vDash^\chi \phi $.
\end{lemma}
\begin{proof}
	Since  $\langle H^\chi \rangle^\chi = H^\chi $  and by the fact that $\langle H^\chi \rangle$ is a subalgebra of $H$, it follows that $\llbracket\psi \rrbracket^{\langle H^\chi\rangle}_{\sigma}=\llbracket\psi \rrbracket^H_{\sigma}$, for any formula $\psi$, from which the result follows.
\end{proof}

\begin{proposition}\label{lem10}
	Let $H$ be a Heyting algebra and $L$ an intermediate logic. Then we have that $H\vDash^\chi L^\chi$ entails $\langle H^\chi \rangle\vDash L$.
\end{proposition}
\begin{proof}

	Consider any Heyting algebra $H$, and suppose  that $\langle H^\chi \rangle\nvDash L$, then there is some formula $\phi\in L$ and some valuation $V$ such that $(\langle H^\chi \rangle,V)\nvDash\phi$.
	Now, since $\langle H^\chi \rangle$ is the subalgebra generated by $H^\chi$, we can express every element $x\in \langle H^\chi \rangle$ as a polynomial $\bm{\delta}^x$ of elements of $H^\chi$.
	We thus have $ x= \bm{\delta}^x(\overline{y})   $, where for each $y_i$ we have that $y_i \in H^\chi$.
	By writing $\overline{p}=p_1,..., p_n$ for the variables contained in $\phi$ and $\bm{\delta}(\overline{y})$ for the sequence of polynomials corresponding to the elements  $x_1=V(p_1),..., x_n=V(p_n)$, we get that $ \llbracket \phi(\overline{p}) \rrbracket^{\langle H^\chi \rangle}_V = \bm{\phi}(\bm{\delta}(\overline{y})) $.
	Since all the elements $\overline{y}$ in the polynomials $\bm{\delta}^x$ are fixed points of $\chi$, we can define a $\chi$-valuation $\sigma: \at \rightarrow H^\chi$ such that  $\sigma: q_i \mapsto y_i$ for all $i\leq n$.
	Then it follows immediately that  $ \llbracket \phi[\sfrac{\delta(\overline{q})}{\overline{p}}]   \rrbracket^{\langle H^\chi \rangle}_{\sigma} = \bm{\phi}(\bm{\delta}(\overline{y})) $.
	But then, since we also had $\llbracket \phi(\overline{p})   \rrbracket^{\langle H^\chi \rangle}_V =  \bm{\phi}(\bm{\delta}(\overline{y})) $, it follows that $ \llbracket \phi[\sfrac{\delta(\overline{q})}{\overline{p}}]   \rrbracket^{\langle H^\chi \rangle}_{\sigma}= \llbracket \phi(\overline{p})   \rrbracket^{\langle H^\chi \rangle}_V $.
	So since $(\langle H^\chi \rangle,V)\nvDash\phi$, we also get that $(\langle H^\chi \rangle,\sigma)\nvDash^\chi\phi[\sfrac{\delta(\overline{y})}{\overline{p}}]$.
	So it then follows by Lemma \ref{lem9.5} that  $H\nvDash^\chi\phi[\sfrac{\delta(\overline{q})}{\overline{p}}] $.
	Now, since $L$ is an intermediate logic, it admits free substitution and so, since $\phi\in L$, we also get that $\phi[\sfrac{\delta(\overline{q})}{\overline{p}}] \in L\subseteq L^\chi$.
	Finally, this means that $ H \nvDash^\chi L^\chi$, thus proving our claim.
\end{proof}

Finally, the former results motivate the introduction of suitable $\chi$-varieties, which we will show being the correct semantic counterpart to $\chi$-logics. Let $\vari$ be an arbitrary variety of Heyting algebras, then its $\chi$-closure is the class:
\begin{equation*}
	\vari^\chi= \{\; K \in \ha \;|\;  \exists H.\; H\in \vari \text{ and } H^\chi=K^\chi  \;\}.
\end{equation*}

\noindent
We say that a Heyting algebra $K$ is a \textit{core superalgebra} of $H$ if $ H^\chi=K^\chi $  and $ H\preceq K  $. We say that $\dvari$ is $\chi$-variety if $\dvari= \vari^\chi$ for some variety $\vari$ of Heyting algebras.
We then prove the following result which characterises $\chi$-varieties. 

 \begin{theorem}\label{nvtv}
	A class of Heyting algebras $\mathcal{C}$  is a $\chi$-variety if and only if it is closed under subalgebras, homomorphic images, products and core superalgebras.
\end{theorem}
\begin{proof}
	($\Leftarrow$) Suppose $\mathcal{C}$ is closed under subalgebras, homomorphic images, products and core superalgebras. Obviously $\mathcal{C}$ is a variety and for any Heyting algebra $H$ such that there is some $K\in \class{C}$ with $H^\chi=K^\chi$ and $K\preceq H$, it follows by closure under core superalgebra that $H\in\class{C}$. Therefore, it follows that $\mathcal{C}=\mathcal{C}^\chi$, hence $\class{C}$ is a $\chi$-variety. ($\Rightarrow$) Suppose $\mathcal{C}$  is a $\chi$-variety, i.e. $\class{C}=\nvari$ for some variety $\vari$. We show that $\class{C}$ is closed under subalgebras, as the other cases follow by an analogous reasoning. Suppose $H\in \class{C}$ and $K\preceq H$. Since $\class{C}=\nvari$ there is some $H'\in \vari$ such that $(H')^\chi=H^\chi$ and $H'\preceq H$. Then consider $K'=K\cap H'$. Since the intersection of two subalgebras is still a subalgebra and since $K'\subseteq H'$, it follows that $K'\preceq H'$ and therefore $K'\in \vari$. Moreover, by a similar reasoning we have that $K'\preceq K$. Finally, since $(K')^\chi = K^\chi \cap (H')^\chi $ and $(H')^\chi = H^\chi$, we have $(K')^\chi = K^\chi \cap H^\chi = K^\chi $. Therefore, by the fact that $K'\preceq K$, $(K')^\chi =  K^\chi $ and $K\in \vari$, we obtain that $K\in\nvari=\class{C}$.	
\end{proof}

\noindent
It is then easy to show that $\chi$-varieties form a bounded lattice with operations $\dvari_0\land \dvari_1:=\dvari_0\cap \dvari_1$ and  $\dvari_0\lor \dvari_1:=\bm{\dvari}(\dvari_0\cup \dvari_1)$, where $\bm{\dvari}(\mathcal{C})$ denotes the smallest $\chi$-variety containing $\mathcal{C}$.
We shall denote the lattice of $\chi$-varieties by $\HA^{\chi}$.
One can show that the map $\vari\mapsto\vari^\chi$ is a lattice homomorphism.
Together with the results of the previous sections, we have thus obtained a lattice $\IL^{\chi}$ of $\chi$-variants of intermediate logics, and a lattice $\HA^{\chi}$ of $\chi$-varieties.
In the next section we shall see how to relate these two structures in order to prove the completeness of the algebraic semantics we introduced.

\section{Duality}\label{dual}

In this section we shall show that the lattice of $\chi$-logics $\IL^{\chi}$ and the lattice of $\chi$-varieties $\HA^{\chi}$ are dual to each other. We prove this result by relying on the standard dual isomorphism between the lattice of intermediate logics and the lattice of varieties of Heyting algebras. We derive as corollaries of such isomorphism a completeness theorem for $\chi$-logics. 

Let $\Gamma$ be a set of formulas and $\class{C}$ a class of Heyting algebras, then we define the two maps $\Var^\chi$ and $\Log^\chi$ such that:
\begin{align*}
&\Var^\chi: \Gamma \mapsto \{\; H\in \HA  \;|\;  H\vDash^{\chi} \Gamma  \;\};  \\	
&\Log^\chi: \class{C} \mapsto \{\; \phi \in \formulas \;|\;   \class{C}\vDash^{\chi} \phi  \;\}.
\end{align*} 

\noindent A class of Heyting algebras $ \mathcal{C} $ is \emph{$\chi$-definable} if there is a set $\Gamma$ of formulas such that $\mathcal{C} = \Var^\chi(\Gamma)$. We say that a $\chi$-logic $\Lambda$ is \emph{algebraically complete} with respect to a class of Heyting algebras $\mathcal{C}$ if $\Lambda = \Log^\chi(\mathcal{C})$. We will next show that $\Var^\chi(\Gamma)$ is always a $\chi$-variety and $\Log^\chi(\mathcal{C})$ is always a $\chi$-logic. This will later allow us to consider $\Log^\chi$ and $\Var^\chi$ as maps between the lattices  $\IL^{\chi}$ and  $\HA^\chi$.

\begin{proposition}\label{proposition:pippoBello}
	$\chi$-validity is preserved by taking subalgebras, products, homomorphic images and core superalgebras.
\end{proposition}
\begin{proof}
	\emph{(Subalgebras)} Suppose by contraposition $K\preceq H$ and $(K,\sigma)\nvDash^\chi\phi$ for some $\chi$-valuation $\sigma$, then obviously $(H,\sigma)\nvDash^\chi\phi$.
	\emph{(Products)} Let $f:H\twoheadrightarrow K$ be a surjective morphism. If $K\nvDash^\chi\phi$, then by Proposition \ref{lem8.5} it follows that $K\nvDash \phi[\sfrac{\chi^n(\overline{p})}{\overline{p}}]$. Since validity is preserved by homomorphic images, we have that $H\nvDash \phi[\sfrac{\chi^n(\overline{p})}{\overline{p}}] $ and therefore, by Proposition \ref{lem8.5}, $H\nvDash^\chi\phi$.
	\emph{(Homomorphic image)} For products we need to show that if $ \prod_{i\in I} A_i \vDash^\chi \phi$, then $A_i\vDash^\chi \phi$ for all $i\in I$. This claim follows immediately by noticing $\left( \prod_{i\in I} A_i \right)_{\chi} = \prod_{i\in I} \left( A_i \right)_{\chi}$, and so $\chi$-valuations over $\prod_{i\in I} A_i$ are all and only the function-products of $\chi$-valuations over the $A_i$.
	\emph{(Core superalgebra)} Let  $K^\chi =H^\chi$ and $H\preceq K$. By \emph{reductio} suppose that $K\nvDash^\chi\phi$. Then for some valuation $\sigma$ we have  $(K,\sigma)\nvDash^\chi\phi$. Since $H^\chi=K^\chi$ and $H\preceq K$, $\sigma$ is a valuation over $H$ and $\llbracket\phi\rrbracket^H_{\sigma}=\llbracket\phi\rrbracket^K_{\sigma}\neq 1$.
\end{proof}

\begin{corollary}\label{nv1Cor}
	For every set of formulas $\Gamma$, the class of Heyting algebras $\Var^\chi(\Gamma)$ is a $\chi$-variety.
\end{corollary}
\begin{proof}
	It follows from Theorem \ref{nvtv} and Proposition \ref{proposition:pippoBello}.
\end{proof}

\noindent It is a straightforward consequence of Corollary \ref{nv1Cor} that every $\chi$-definable class of Heyting algebras is also a $\chi$-variety.
The next proposition shows that for every class $\class{C}$ of Heyting algebras its set of validities $\Log^\chi(\class{C})$ is a $\chi$-logic. 

\begin{proposition}\label{nl}
	For every set of algebras $\mathcal{C}$, the class of formulas $\Log^\chi(\mathcal{C}) $ is a $\chi$-logic.
	Moreover, $\Log^\chi(\mathcal{C})$ is the $\chi$-variant of $\Log(\mathcal{C})$
\end{proposition}
\begin{proof}
	We have:
	\begin{align*}
	\phi \notin \Log^\chi(\class{C}) & \Longleftrightarrow \exists H\in \class{C} \text{ such that } H\nvDash^\chi \phi & \\
	& \Longleftrightarrow \exists H\in \class{C} \text{ such that } H\nvDash\phi[\sfrac{\chi^n(\overline{p})}{\overline{p}}] & \text{(by Proposition $\ref{lem8.5}$)}\\
	& \Longleftrightarrow \phi[\sfrac{\chi^n(\overline{p})}{\overline{p}}]\notin \Log(\class{C}) &\\
	& \Longleftrightarrow \phi \notin (\Log(\class{C}))^\chi. &
	\end{align*}
	\noindent Hence  $\Log^\chi(\class{C}) $ is the $\chi$-variant of $\Log(\class{C})$.
\end{proof}

\noindent
These two results establish that the maps $\Log^\chi: \HA^{\chi} \rightarrow \IL^{\chi}$ and $\Var^\chi: \IL^{\chi} \rightarrow \HA^{\chi}$ are well-defined.
Using the homomorphisms $\vari\mapsto \vari^\chi$ and $L \mapsto L^\chi$, and the isomorphism $\HA \cong^{op} \IL$ between intermediate logics and varieties of Heyting algebras given by the maps $\Log$ and $\Var$, we obtain the following commuting diagrams.

\begin{proposition}\label{prop3.6}
	For every intermediate logic $L$, $\Var^{\chi}(L^{\chi})=  \Var(L)^\chi$.
\end{proposition}
\begin{center}
	\begin{tikzcd}
		\IL \arrow[dd, "\Var"'] \arrow[rrr] \arrow[rrr, "(-)^\chi"] &  &  & \IL^{\chi} \arrow[dd, "\Var^\chi"] \\
		&  &  &                           \\
		\HA \arrow[rrr, "(-)^\chi"]                            &  &  & \HA^{\chi}                       
	\end{tikzcd}
\end{center}
\begin{proof}
	$(\subseteq)$ Consider any Heyting algebra $H\in \Var^{\chi}(L^{\chi})$. Then we have $H\vDash^\chi L^{\chi}$ and by Proposition \ref{lem10} it follows  $\langle H^\chi\rangle\vDash  L  $. So we clearly have that $\langle H^\chi\rangle\in \Var(L)$ and since $\langle H^\chi\rangle^\chi=H^\chi$ and $\langle H^\chi\rangle\preceq H$ also $H\in \Var(L)^\chi$. 	$(\supseteq)$ Consider any Heyting algebra $H\in \Var(L)^\chi$, then there is some $K\in \Var(L)$ such that $K\preceq H$ and $H^\chi=K^\chi$. Then we have that $K\vDash L$, so by Corollary \ref{lem9} above $K\vDash^\chi L^\chi$ which entails $K\in \Var^\chi(L^\chi)$. Finally, since $\chi$-varieties are closed under core superalgebra, it follows that $H\in \Var^\chi(L^\chi)$.
\end{proof}

\begin{proposition}\label{prop3.7}
	For every variety $\vari$ of Heyting algebras $\Log^{\chi}(\vari^\chi)= \Log(\vari)^\chi$.
\end{proposition}
\begin{center}
	\begin{tikzcd}
		\IL \arrow[rrr] \arrow[rrr, "(-)^\chi"]           &  &  & \IL^{\chi}                            \\
		&  &  &                            \\
		\HA \arrow[rrr, "(-)^\chi"] \arrow[uu, "\Log"] &  &  & \HA^{\chi}     \arrow[uu, "\Log^\chi"']
	\end{tikzcd}
\end{center}
\begin{proof}
We prove both directions by contraposition. $(\subseteq)$ Suppose $\phi\notin \Log(\vari)^\chi$, then  $\phi[\sfrac{\chi^n(\overline{p})}{\overline{p}}]\notin \Log(\vari)$ and hence there is some $H\in \vari$ such that  $H \nvDash \phi[\sfrac{\chi^n(\overline{p})}{\overline{p}}]$. By Proposition \ref{lem8.5}  $H \nvDash^\chi \phi$, hence $\phi \notin \Log^{\chi}(\vari^\chi)$. 	$(\supseteq)$ Suppose $\phi\notin \Log^{\chi}(\nvari)$. It follows that there is some $H\in \nvari$ such that $H\nvDash^\chi \phi$, hence by Lemma \ref{lem9.5}  $\langle H^\chi \rangle\nvDash^\chi \phi$. It thus follows by Proposition \ref{lem8.5} that $\langle H^\chi \rangle\nvDash \phi[\sfrac{\chi^n(\overline{p})}{\overline{p}}]$. 	Now, since  $H\in \nvari$, we have for some $K\in\vari$ that $K\preceq H$ and $K^\chi= H^\chi$. Thus it follows that $\langle H^\chi \rangle\preceq K$ and therefore $\langle H^\chi\rangle\in \vari$. Finally, since $\langle H^\chi \rangle\nvDash \phi[\sfrac{\chi^n(\overline{p})}{\overline{p}}]$ we get that $\phi[\sfrac{\chi^n(\overline{p})}{\overline{p}}]\notin \Log(\vari)$ and hence $\phi\notin \Log(\vari)^\chi$.	
\end{proof}

\noindent
Since the diagrams above commute, it is easy to prove a definability theorem and a completeness theorem for $\chi$-logics and $\chi$-varieties.

\begin{theorem}[Definability Theorem]\label{DT}
	$\chi$-varieties are defined by their $\chi$-validities: $H\in \dvari$ if and only if $H\vDash^\chi \Log^\chi(\dvari)$.
\end{theorem}
\begin{proof}
	For any  $\chi$-variety $\dvari$ such that $\dvari=\nvari$  we have:
	\begin{align*}
	\Var^\chi(\Log^\chi (\nvari)) &=  \Var^\chi(\Log (\vari)^\chi) &  \text{(by Proposition \ref{prop3.7})}\\	
	&= \Var(\Log (\vari))^\chi  & \text{(by Proposition \ref{prop3.6})}\\
	&=\vari^\chi  & \text{(by standard duality)}		
	\end{align*}
	\noindent Hence $\Var^\chi\circ \Log^\chi = \mathbb{1}_{\HA^{\chi}}$, which proves our claim.
\end{proof}

\begin{theorem}[Algebraic Completeness]\label{ACDNAV}
	$\chi$-logics are complete with respect to their corresponding $\chi$-variety: $\phi\in \Lambda$ if and only if $\Var^\chi(\Lambda)\vDash^\chi \phi$.
\end{theorem}
\begin{proof}
	For any  $\chi$-logic $\Lambda$ such that $\Lambda=L^\chi$ we have:
	\begin{align*}
	\Log^\chi(\Var^\chi(L^\chi)) &= \Log^\chi( \Var(L)^\chi)  & \text{(by Proposition \ref{prop3.6})}\\
	&=\Log(\Var(L))^\chi  & \text{(by Proposition \ref{prop3.7})}\\
	&=L^\chi  & \text{(by standard duality)}
	\end{align*}
	\noindent Hence $\Log^\chi \circ \Var^\chi = \mathbb{1}_{\IL^{\chi}}$, which proves our claim.
\end{proof}

\noindent The former completeness theorem shows that the algebraic semantics that we have introduced in the previous section is indeed a suitable framework to study $\chi$-variants of intermediate logics from a semantics point of view.
Similarly, the definability theorem for $\chi$-varieties allows us to give a first  \textit{external} characterisation of $\chi$-varieties: they are exactly the classes of Heyting algebras which are $\chi$-definable.
Finally, since $\Var^\chi$ and $\Log^\chi$ are lattice homomorphisms, we obtain the following dual isomorphism result.

\begin{theorem}[Duality]\label{DUAL}
	The lattice of $\chi$-logics is dually isomorphic to the lattice of $\chi$-varieties of Heyting algebras, i.e. $\IL^{\chi} \cong^{op} \HA^{\chi}$.
\end{theorem}

Before turning to the study of specific $\chi$-logics in the next section, we shall also provide an alternative characterisation of $\chi$-varieties.
First notice that, since $\chi$-varieties are closed under subalgebras, homomorphic images and products, they are also (standard) varieties and thus Birkhoff Theorem tells us that they are generated by their collection of subdirectly irreducible elements.
It is possible to show more and give an \emph{internal} characterisation of $\chi$-varieties: they are exactly the classes of Heyting algebras generated (also under the core superalgebra operation) by their collection of core generated, subdirectly irreducible elements.
Recall that, given $\class{C}$ a class of Heyting algebras, we indicate with $\bm{\dvari}(\class{C})$ the least $\chi$-variety containing $\class{C}$ and with $\bm{\vari}(\class{C})$ the least variety containing $\class{C}$.
We first adapt Tarski's HSP-Theorem to the current setting. 

\begin{theorem}\label{TNV}
	Let  $\class{C}$ be a class of Heyting algebras, then  $\bm{\dvari}(\class{C})=(\,\Hvar\Svar\Pvar(\class{C})\,)^\chi$.
\end{theorem}
\begin{proof}
	By definition we have $\bm{\dvari}(\class{C})=\bm{\vari}(\class{C})^\chi$ and by Tarski's HSP-Theorem  	$\bm{\vari}(\class{C})=\Hvar\Svar\Pvar(\class{C})$. It immediately follow $\bm{\dvari}(\class{C})=(\,\Hvar\Svar\Pvar(\class{C})\,)^\chi$.
\end{proof}
 
\noindent From the former theorem it is easy to prove the following useful result.

\begin{proposition}\label{31}
	Let $\dvari$ be a $\chi$-variety, then $\dvari=\bm{\dvari}(\class{C}) $ iff $ \Log^\chi(\dvari)=\Log^\chi (\class{C}).$
\end{proposition}
\begin{proof}
	$(\Rightarrow)$ Since $\class{C}\subseteq \dvari$, the inclusion from right to left is straightforward. Suppose now that $\dvari\nvDash^\chi \phi  $ then there is some $H\in\dvari$ such that $H\nvDash^\chi \phi$.
	Then since $\dvari=\bm{\dvari}(\class{C})$, it follow by Theorem \ref{TNV}  that $H\in \Hvar\Svar\Pvar(\class{C})^\chi$.
	By Proposition \ref{proposition:pippoBello}, it follows that for some $A \in \class{C}$ we have $A\nvDash^\chi \phi$.
	Hence $\phi\notin \Log^\chi (\class{C})$.
	
	$(\Leftarrow)$ Suppose  $\Log^\chi(\dvari)=\Log^\chi(\class{C})$.
	It follows that $\Var^\chi(\Log^\chi(\dvari))=\Var^\chi(\Log^\chi(\class{C}))$, hence by the Duality Theorem \ref{DUAL}, we have $\dvari=\Var^\chi(\Log^\chi(\class{C}))$.
	Finally, since $\Log^\chi(\class{C}) = \Log^\chi(\bm{\dvari}(\class{C}))$ by Proposition \ref{proposition:pippoBello} and Theorem \ref{TNV}, we have $\Var^\chi(\Log^\chi(\class{C}))=\Var^\chi(\Log^\chi(\bm{\dvari}(\class{C}) )$;
	and by Duality $\Var^\chi(\Log^\chi(\bm{\dvari}(\class{C}) )=\bm{\dvari}(\class{C}) $, it follows $\dvari=\bm{\dvari}(\class{C}) $.
\end{proof}

\noindent  A first characterisation is given by following result, stating that every $\chi$-variety $\dvari$ is generated by its collection of core generated Heyting algebras. We denote by $\dvari_{CG}$  the subclass of core generated Heyting algebras of a $\chi$-variety $\dvari$.

\begin{proposition}\label{reg}
	Every $\chi$-variety is generated by its collection of core generated elements, i.e. $\dvari = \bm{\dvari}(\dvari_{CG}) $.
\end{proposition}
\begin{proof}
	Let $\dvari$ be a $\chi$-variety, then for any non core generated $H\in \dvari $  we have  $\langle H^\chi \rangle \preceq H$ and $H^\chi=\langle H^\chi \rangle^\chi$. So since $\langle H^\chi \rangle \in \dvari_{CG}$  it follows $H\in \bm{\dvari}(\dvari_{CG}) $.
\end{proof}

\noindent By Birkhoff theorem we know that every $\chi$-variety is generated by its subdirectly irreducible elements and, by the previous proposition, we know that every $\chi$-variety is generated by its core generated elements. The following result shows that the intersection of these two classes of generators actually suffices. If $\dvari$ is a $\chi$-variety, let $\dvari_{CGSI}$  its subclass of core generated subdirectly irreducible Heyting algebras. We prove a version of Birkhoff Theorem for $\chi$-varieties showing that $\dvari=\bm{\dvari}(\dvari_{CGSI})$.

\begin{theorem}\label{birkhoff}
	Every $\chi$-variety is generated by its collection of core generated subdirectly irreducible elements: $\dvari = \bm{\dvari}(\dvari_{CGSI}) $.
\end{theorem}
\begin{proof}
	By the dual isomorphism between $\chi$-logics and $\chi$-varieties it suffices to show that $\Log^\chi (\dvari)=\Log^\chi(\bm{\mathcal{X}}(\dvari_{CGSI})) $.  By Proposition \ref{31} this is equivalent to $\Log^\chi (\dvari)=\Log^\chi(\dvari_{CGSI}) $. The direction  $\Log^\chi (\dvari)\subseteq \Log^\chi(\dvari_{CGSI}) $  follows immediately from the inclusion $\dvari_{CGSI}\subseteq \dvari$. We next show that $\Log^\chi(\dvari_{CGSI}) \subseteq \Log^\chi (\dvari)$. 
	
	Suppose by contraposition $\phi\notin \Log^\chi (\dvari) $, then for some $H\in \dvari$ and some $\chi$-valuation $\sigma$, we have $(H,\sigma)\nvDash^\chi \phi$ and so by Lemma \ref{lem9.5}  $(\langle H^\chi \rangle, \sigma) \nvDash^\chi \phi$. Now, it is a well-known fact, originally shown by Wronski in \cite{Wronski}, that for any Heyting algebra $B$ and $x\in B$  such that $b\neq 1_B$, there is a subdirectly irreducible algebra $C$ and a surjective homomorphism $h:B\twoheadrightarrow C$ such that $f(b)=s_C$, where $s_C$ is the second greatest element of $C$. Then, since $x=\llbracket \phi \rrbracket^{\langle H^\chi\rangle} _{\sigma}\neq 1_H$  there is a subdirectly irreducible algebra $C$ and surjective homomorphism $h:\langle H^\chi \rangle\twoheadrightarrow C$ with $h(x )=s_C$. Consider now the valuation $\tau = h\circ \sigma $ then, since $h$ a is homomorphism, $\tau$ is still a $\chi$-valuation. Let $p_0,\dots,p_n$ be the variables in $\phi$, it  follows  by the properties of homomorphisms that:
	\begin{align*}
	\llbracket \phi(p_0,\dots,p_n) \rrbracket^C _{\tau	} &= \phi_C [\tau(p_0),\dots,\tau(p_n)] \\
	&= \phi_C[ h(\sigma(p_0)),\dots,h(\sigma(p_n)) ] \\ &=  h \llbracket \phi(p_0,\dots,p_n) \rrbracket  ^{\langle H^\chi\rangle} _{\sigma}\\ & = s_C.
	\end{align*}
	
	\noindent From which it  immediately follows that $ (C,\tau)\nvDash \phi$ and so that $C\nvDash\phi$. Now, since $H\in\dvari$, we have that $\langle H^\chi\rangle\in \dvari$ and so since $h:\langle H^\chi\rangle\twoheadrightarrow C$ also that $C\in \dvari$. Moreover, we have that  $C$ is subdirectly irreducible and, since $C=h[\langle H^\chi\rangle]$, also that $C$ is core generated. Finally, this means that $C\in \dvari_{CGSI} $ and so that $\phi \notin \Log^\chi(\dvari_{CGSI}) $, which proves our claim.
\end{proof}

\section{The lattices of $\chi$-logics}\label{discu}

In this Section we will consider examples of $\chi$-logics and look at their specific properties and characterisation.
Recall that in Lemma \ref{lemma:finiteFixpoints} we have shown that there are only $6$ fix-points of intuitionistic univariate formulas: $\bot, p, \neg p, \neg \neg p, p\lor \neg p, \top$.
Since $L^\neg = L ^{\neg \neg}$ for every intermediate logic $L$, this means that there are at most five lattices of $\chi$-logics.
We provide a description of these lattices.

\smallskip\noindent
\textbf{$p$-logics:}
Firstly, the lattice of $p$-logics $ \IL^{p} $ actually coincides with the lattice of intermediate logics $ \IL$, since for every intermediate logic $L$ it is clearly the case that $L^p=L$.
From the algebraic perspective, this means that for any Heyting algebra $H$ its $p$-core is $H^p=H$, thus we are not imposing any restriction on our valuations.

\smallskip\noindent
\textbf{$\top$-logics and $\bot$-logics:}
The two ``limit'' cases  $ \IL^{\bot} $ and $\IL^{\top}$ are more interesting.
Notice that $H^\bot = \{ 0_H \}$ and $H^{\top} = \{ 1_H \}$, and so under the algebraic semantics that we have introduced $\bot$-models allow only the constant valuation with image $0_H$ and, analogously, $\top$-models allow only the constant valuation with image $1_H$.
Interestingly, this means that the notion of core superalgebra collapses in both cases to that of superalgebra, as we have $\langle H^{\bot} \rangle  =  \langle H^{\top} \rangle  =  \{ 0_H, 1_H \}$, which is a subalgebra of every Heyting algebra.

Thus there is only one $\bot$-variety and only one $\top$-variety, in both cases the variety of all Heyting algebras.
By the duality result of the previous section, this means there are exactly one $\bot$-logic ($\IPC^{\bot}$) and one $\top$-logic ($\IPC^{\top}$), which are respectively the $\bot$-variant and $\top$-variant of every intermediate logic.
These two logics are characterised by the following properties:
\begin{equation*}
\begin{array}{l @{\hspace{2em}\text{iff}\hspace{2em}} l @{\hspace{2em}\text{iff}\hspace{2em}} l}
	\phi(p_1,\dots,p_n) \in \IPC^{\bot}
		&  \phi(\bot,\dots,\bot) \in \IPC
		&  \phi(\bot,\dots,\bot) \in \CPC \\
	\phi(p_1,\dots,p_n) \in \IPC^{\top}
		&  \phi(\top,\dots,\top) \in \IPC
		&  \phi(\top,\dots,\top) \in \CPC
\end{array}
\end{equation*}
Notice in particular that, although they correspond to the same variety, the two logics are distinct.

\smallskip\noindent
\textbf{$\neg p$-logics:}
Apart from $\IL$, the lattice  $ \IL^{\neg p} $ is the only lattice of $\chi$-logics that has already been studied in the literature, although under a different name.
In fact an example of $\neg p$-logic is \emph{inquisitive logic} $\inqB$, which is the $\neg p$-variant of the intermediate logics $\KP,\ND $ and $\ML$ as shown in Theorem \ref{inq}.
As a matter of facts, the algebraic semantics for inquisitive logic restricting valuations of atomic formulas to regular elements was already introduced in \cite{BezhanishviliGrillettiHolliday2019}, and was later generalised in \cite{Quadrellaro.2019} to consider the entire lattice of $\mathtt{DNA}$-logics\footnote{Also referred to as \emph{negative variants} in the literature (\cite{Miglioli:89,Iemhoff:2016aa,Ciardelli:09thesis}).} and their corresponding varieties.
This semantics coincides with the one introduced in this paper, since regular elements are exactly the fixpoints of the $\neg\neg$ operator.
This approach has proved to be particularly useful:
for instance, \cite{Quadrellaro.2019} shows that the lattice of extensions of $\inqB$ is dually isomorphic to $\omega + 1$, and also provide an axiomatisation of all such extensions by a generalisation of the method of Jankov formulas.

$\neg p$-logics have a particularly interesting feature:
as mentioned before, the $\neg p$-core of a Heyting algebra is the set of its regular elements, which is a Boolean algebra for the signature $\{1,0,\land,\rightarrow\}$.
This easily entails the following corollary: Given an intermediate logic $L$ and a $\vee$-free formula $\phi$, $\phi \in L^{\neg p}$ iff $\phi$ is a classical tautology (Theorem 2.5.2 in \cite{Ciardelli:16}).
That is, $\neg p$-logics are are logics whose $\{1,0,\land,\rightarrow\}$-fragment behaves classically, and which present an intuitionistic behaviour once formulas containing disjunction are concerned.
Such intuitionistic behaviour disappears once also disjunction is forced to be classical, as the following lemma shows:
\begin{proposition}\label{prop:negativeWemCpc}
	Let $L$ be an intermediate logic.
	Then $L^{\neg p } = \CPC$ iff $L$ extends the logic of week excluded middle $\WEM := \IPC + (\neg p \vee \neg\neg p)$.
\end{proposition}

\noindent
The original proof of this result is given in \cite{Ciardelli:09thesis} (Proposition 5.2.22).
Here we present an alternative proof using the machinery developed in the previous sections.

\begin{proof}
	Firstly, notice that $L^{\neg p} = \CPC$ iff $q\vee \neg q \in L^{\neg p}$.
	The left-to-right implication is trivial.
	As for the other implication, $q \vee \neg q \in L^{\neg p}$ implies that the set of regular elements of an algebra in the variety $\bm{\mathcal{X}}(L^{\neg p})$ is itself a subalgebra, and a Boolean algebra at that.
	By Proposition \ref{reg} and Theorem \ref{ACDNAV}, it easily follows that $L^{\neg p} = \CPC$.

	The main statement now follows easily: By Proposition \ref{lemma:characterizationChiVariant}, $q \vee \neg q \in L^{\neg q}$ iff $\neg q \vee \neg\neg q \in L$, which in turn is equivalent to $\mathtt{WEM} \subseteq L$.
\end{proof}

\noindent
We refer the reader to \cite{Quadrellaro.2019} for more information on $\neg p$-logics and $\neg p$-varieties.

\smallskip\noindent
\textbf{$(p\vee \neg p)$-logics:}
Finally, let us consider the lattice $ \IL^{p\lor \neg p} $. The next proposition gives a characterisation of the $p\lor \neg p$-core of any Heyting algebra $H$.
 
 \begin{proposition}
 	Let $H$ be a Heyting algebra and let $x\in H$. The following are equivalent:
 	\begin{enumerate}
 		\item $x = y \lor \neg y$ for some $y\in H$;
 		\item $\neg x = 0$;
 		\item for every $y\in H$, if $x\land y = 0$, then $y=0$.
 	\end{enumerate}
 \end{proposition}
 \begin{proof}
 	$(1\Rightarrow 2)$ Suppose $x = y \lor \neg y$ for some $y\in H$.
 	Then $\neg x = \neg (y \lor \neg y)=\neg y \land \neg \neg y = 0_H$.
 	$(2\Rightarrow 3)$  Suppose $\neg x = 0$ and $x\land y = 0$.
 	Then $\neg \neg x \land \neg \neg y = \neg \neg 0$, hence $1 \land \neg \neg y = \neg\neg y = 0$.
 	Since $y\leq \neg\neg y$, it follows that $y=0$.
 	$(3\Rightarrow 1)$ Suppose $x$ is as in point $3$.
 	Since $x\land \neg x = 0$, it follows that $\neg x = 0$. Consequently, we also have $x=x\lor \neg x$.
 \end{proof}
 

\begin{wrapfigure}{r}{0.29\textwidth}
	\centering
	\begin{tikzpicture}
		\node[highlight] at ( 0,0) {};
		\node[highlight] at ( 0,2) {};
		\node[highlight] at ( 0,3) {};

		\node[dot,label={0:$0$}] (0) at ( 0,0) {};
		\node[dot,label={180:$a$}] (a) at (-1,1) {};
		\node[dot,label={0:$b$}] (b) at ( 1,1) {};
		\node[dot,label={0:$s$}] (s) at (0,2) {};
		\node[dot,label={0:$1$}] (1) at (0,3) {};

		\draw (0) -- (a);
		\draw (0) -- (b);
		\draw (a) -- (s);
		\draw (b) -- (s);
		\draw (s) -- (1);
	\end{tikzpicture}
	\caption{An example of an algebra in $\Var^{p\vee \neg p}(\LC^{p\vee\neg p})$ but not in $\Var(\LC^{p\vee\neg p})$. The circles indicate the members of the subalgebra generated by the dense elements. Notice that this algebra is not core-generated.}
	\label{figure:LC}
\end{wrapfigure}
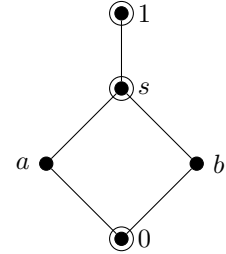

\noindent
The elements satisfying properties 1,2 and 3 above are referred to as \textit{dense elements}.
Notice that property 1 is exactly the condition defining the elements of $H^{p\lor\neg p}$, thus the previous proposition provides a characterisation of the $(p\lor\neg p)$-core of $H$.

Now, it is easy to see that the dense elements of a Heyting algebra form a \textit{filter} and that they are closed under the operations $\land, \lor, \rightarrow$ and $1$.
As a simple consequence of this, we have that for any Heyting algebra $H$ its core subalgebra is  $\langle H^{p\lor\neg p}\rangle= H^{p\lor\neg p} \cup \{0 \}$.
Therefore, the core generated algebras---which by Theorem \ref{birkhoff} suffice to generate all the $(p \lor\neg p)$-varieties---are exactly the algebras containing only dense elements apart from $0$.

 

We obtain an interesting example of $(p\lor \neg p)$-logic by taking the $(p\lor\neg p)$-variant of \emph{G\"{o}del-Dummett logic} $\LC$.
Recall that $\LC$ is the intermediate logic extending $\IPC$ with the axiom $(p\rightarrow q) \lor (q\rightarrow p)$.
It can be also characterised as the logic of linear Heyting algebras (Example 4.15 in \cite{Chagrov:97}).
Analogously, the logic $ \LC^{p\lor \neg p} $ forces a similar linearity condition, but now limited to the dense elements of a Heyting algebras.
Notice that by Proposition \ref{nl}, the variety $\Var^{p\lor\neg p}(\LC^{p\lor \neg p})$ is still generated by the class of linear Heyting algebras---however, the closure under core-superalgebras leads to a variety properly extending $\Var(\LC)$, as shown in Figure \ref{figure:LC}.
Moreover, notice that linear algebras are core-generated since $\neg x = 0$ for every non-zero element $x$;
thus we found a class of core-generated algebras which generate the whole $(p\vee\neg p)$-variety.

Finally, we have seen in Proposition \ref{prop:negativeWemCpc} that the intermediate logics whose $\neg p$-variant is $\CPC$ are exactly the extensions of $\mathtt{WEM}$.
So a natural question to ask is what intermediate logics have $\CPC$ as their $(p\vee\neg p)$-variant.
The next proposition establishes that $\CPC$ itself is the only intermediate logic with this property.

\begin{proposition}
	$L^{p\lor \neg p}= \mathtt{CPC}$ iff $L=\mathtt{CPC}$.
\end{proposition}
\begin{proof}
	The left-to-right implication is trivial.
	As for the other implication, suppose $L^{p\lor \neg p}= \mathtt{CPC}$.
	In particular, $q\lor \neg q \in \CPC = L^{p\lor \neg p}$.
	By Lemma \ref{lemma:characterizationChiVariant}, we have $(q\vee\neg q) \vee \neg( q\vee\neg q) \equiv q\vee\neg q \in L$, which means $L=\CPC$.
 \end{proof}


\smallskip\noindent
Now that we have described the five lattices of $\chi$-logic more in detail, we are ready to show they are distinct.
However, we need to clarify what we mean by distinct lattices:
as we have already seen, $\IL^{\top}$ and $\IL^{\bot}$ both contain only one logic, so these lattices are isomorphic;
but we already observed that the logics $\IPC^{\top}$ and $\IPC^{\bot}$ are different, so we will consider these lattices distinct.
More generally, we will say that $\IL^{\chi}$ and $\IL^{\theta}$ are equal if for every intermediate logic $L$ we have $L^{\chi} = L^{\theta}$---as in the case of $\IL^{\neg p}$ and $\IL^{\neg\neg p}$;
otherwise we will say the lattices are distinct.

This suggests to study the relation between these lattices in a more systematic way:
define the pointwise extension relation $\IL^{\chi} \latleq \IL^{\theta}$ to hold if for every intermediate logic $L$ we have $L^{\chi} \subseteq L^{\theta}$.
The relation $\latleq$ is a partial order between the lattices of $\chi$-logics.
In particular, $\IL^{\chi} \latleq \IL^{\theta} \latleq \IL^{\chi}$ if and only if the two lattices $\IL^{\chi}$ and $\IL^{\theta}$ are equal.
The following Theorem characterises the properties of this relation.

\begin{theorem}\label{theorem:orderLattices}
	Let $\chi$ and $\theta$ be univariate formulas.
	Then the following are equivalent:
	\begin{enumerate}
		\item $\IL^{\chi} \latleq \IL^{\theta}$;
		\item $\IPC^{\chi} \subseteq \IPC^{\theta}$;
		\item $(\theta^{2}(p) \leftrightarrow p) \rightarrow (\chi^{2}(p) \leftrightarrow p) \in \IPC$;
		\item For every Heyting algebra $H$, $H^{\theta} \subseteq H^{\chi}$.
	\end{enumerate}
\end{theorem}

\begin{proof}
	$(1 \Rightarrow 2)$
	It follows from the definition of $\latleq$.
	$(2 \Rightarrow 3)$
	Since $\IPC^{\chi} \subseteq \IPC^{\theta}$, we have in particular that $\chi^2(p) \leftrightarrow p \in \IPC^{\theta}$.
	This means that $\IPC + (\theta^2(p) \leftrightarrow p) \vDash \chi^{2}(p) \leftrightarrow p$;
	and so by the deduction theorem of $\IPC$ we have $(\theta^{2}(p) \leftrightarrow p) \rightarrow (\chi^{2}(p) \leftrightarrow p) \in \IPC$.
	$(3 \Rightarrow 4)$
	We prove the contrapositive of the implication: suppose that $H^{\theta} \nsubseteq H^{\chi}$ for some Heyting algebra $H$.
	Consider an element $a \in H^{\theta} \setminus H^{\chi}$.
	By Lemma \ref{lemma:fixpointsAlgebras} we have $\bm{\theta}^2(a)\leftrightarrow a = 1_H$ and $\bm{\chi}^2(a) \leftrightarrow a \ne 1_H$.
	So in particular $H\nvDash (\theta^{2}(p)\leftrightarrow p) \rightarrow (\chi^2(p)\leftrightarrow p)$, showing that this is not a formula in $\IPC$.
	$(4 \Rightarrow 1)$
	Consider an intermediate logic $L$ and
	take an arbitrary formula $\phi \notin L^{\theta}$.
	By Proposition \ref{nl} $L^{\theta} = \Log^{\theta}(\Var(L))$, and so there exists an algebra $H\in \Var(L)$ and a $\theta$-valuation $\sigma$ such that $(H,\sigma) \nvDash \phi$.
	Since $H^\theta \subseteq H^\chi$ by hypothesis, it follows that $\phi \notin \Log^{\chi}(\Var(L))$ either.
	Again by Proposition \ref{nl}, $L^{\chi} = \Log^{\chi}(\Var(L))$, and thus $\phi \notin L^{\chi}$.
	Since $\phi$ was generic, it follows that $L^{\chi} \subseteq L^{\theta}$, as wanted.
\end{proof}

\bgroup
\begin{wrapfigure}{r}{.33\textwidth}
	\centering
	\begin{tikzpicture}
		\node (p)    at ( 0,0) {$\IL^p = \IL$};
		\node (nnp)  at (-1,1) {$\IL^{\neg\neg p}$};
		\node (ponp) at ( 1,1) {$\IL^{p\vee\neg p}$};
		\node (bot)  at (-1,2) {$\IL^{\bot}$};
		\node (top)  at ( 1,2) {$\IL^{\top}$};

		\draw (p)    -- (nnp);
		\draw (p)    -- (ponp);
		\draw (nnp)  -- (bot);
		\draw (nnp)  -- (top);
		\draw (ponp) -- (top);
	\end{tikzpicture}
	\caption{The Hasse diagram of the 5 lattices of $\chi$-logics, ordered under the relation $\latleq$. The diagram is computed using Theorem \ref{theorem:orderLattices}.}
\end{wrapfigure}
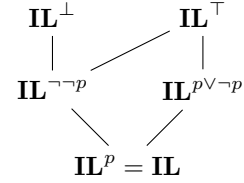

\noindent
\begin{minipage}[t]{.62\textwidth}
\begin{corollary}
	There are exactly 5 lattices of $\chi$-logics, for $\chi$ a univariate formula:
	\begin{equation*}
		\IL^{\bot},\hspace{1.5em}
		\IL^{p} = \IL,\hspace{1.5em}
		\IL^{\neg p} = \IL^{\neg\neg p},\hspace{1.5em}
		\IL^{p \vee \neg p},\hspace{1.5em}
		\IL^{\top}.
	\end{equation*}
\end{corollary}

\begin{proof}
	What remains to be shown is that the lattices are distinct.
	By Theorem \ref{theorem:orderLattices}, we can to this by exhibiting a Heyting algebra $H$ for which the $\chi$-cores are all distinct.
	Indeed, the algebra in Figure \ref{figure:LC} is an example of such an algebra:
	\begin{equation*}
	\begin{array}{ll @{\hspace{3.5em}} ll}
		H^{\bot} &= \{ 0 \}
			&H^{\top} &= \{1\}   \\
		H^{\neg p} &= \{0,a,b,1\}
			&H^{p\vee \neg p} &= \{0,s,1\}   \\
		H^{p} &= H \\
	\end{array}
	\end{equation*}
\end{proof}
\end{minipage}
\egroup

\section{Conclusion}\label{conclusion}

In this article we introduced $\chi$-logics and
a sound and complete algebraic semantics for them, based on Ruitenburg's Theorem.
In Section \ref{chi} we defined the notion of $\chi$-logics and studied them from a syntactical perspective, showing that for a fixed $\chi$ they form a bounded lattice, and that we have only 5 such lattices.
In Section \ref{alge} we defined an algebraic semantics for $\chi$-logics, by relying on an algebraic interpretation of Ruitenburg's Theorem originally described in \cite{santocanale:hal-01766636} and we introduced $\chi$-varieties as the semantic counterpart of $\chi$-logics.
In Section \ref{dual} we proved that $\chi$-logics are indeed complete with respect to their corresponding algebraic semantics and we proved some of their properties.
Finally, in Section \ref{discu}, we have looked more in detail at each of the 5 lattices of $\chi$-logics and characterised explicitly the point-wise extension relation $\latleq$ between the lattices.

The results of this article provide a first approach to generate and study new logics and corresponding algebraic semantics in a systematic fashion.
This work can be extended in several directions:
Firstly, we limited ourselves to univariate formulas but the approach based on Ruitenburg's Theorem can be generalised to the lattice produced by an arbitrary formula---although with a more complex technical machinery.
Secondly, even in this more general setting cores are required to be definable, but it is natural to consider more general notions of core (i.e., fixpoints of some operator) and their corresponding logics.
This step would also allow to move from the setting of intermediate logics to other families of logics.
Another interesting direction of work would be to interpret the results presented in this paper in terms of topological duality, that is, Esakia duality for Heyting algebras (\cite{BezhanishviliEsakia2019}).
Giving a topological interpretation to the results and constructions presented (such as the core-superalgebra operation) would give novel tools to study the structure of the lattices of $\chi$-logics.


\bibliographystyle{abbrv}
\bibliography{tbillc2019}

\appendix
\section{Proof of Lemma \ref{lemma:finiteFixpoints}}\label{appendixA}

\begin{proof}[Proof of Lemma \ref{lemma:finiteFixpoints}]
	As shown in \cite{nishimura_1960,rieger_1949}, the following is a presentation of all the non-constant univariate intuitiornistic formulas modulo logical equivalence.
	\begin{equation*}
		\beta_1 := p
		\hspace{3em}
		\beta_{n+1} := \alpha_n \vee \beta_n
		\hspace{5em}
		\alpha_1 := \neg p
		\hspace{3em}
		\alpha_{n+1} := \alpha_n \to \beta_n
	\end{equation*}

	\noindent
	We will consider the following two properties for a univariate formula $\phi$:
	\begin{enumerate}
		\item $\neg\neg \phi \equiv \top$.
		\item If $\psi$ has property $1$, then $\phi[\sfrac{\psi}{p}] \equiv \top$.
	\end{enumerate}

	\noindent
	In particular, if $\phi$ has both properties then $\phi^2 \equiv \top$, that is, the fix-point of $\phi$ is $\top$.

	Firstly notice that
	\begin{equation*}
		\alpha_5 = ((\neg\neg p) \to p\vee \neg p) \to (\neg p \vee \neg\neg p)
		\qquad\text{and}\qquad
		\beta_5 = ((\neg\neg p) \to p\vee \neg p) \to (\neg\neg p \to p)
	\end{equation*}
	have both properties.

	\vspace{.3em}
	\noindent
	\begin{minipage}[t]{0.48\textwidth}
		$\alpha_5$ has property 1:
		\vspace{-.5em}
		\begin{equation*}
		\begin{array}{rl}
			\neg\neg \alpha_5  \quad
			\scriptstyle\equiv &\scriptstyle\neg\neg( ((\neg\neg p \to p) \to p \vee \neg p) \to (\neg p \vee \neg\neg p) ) \\
			\scriptstyle\equiv &\scriptstyle\neg\neg( (\neg\neg p \to p) \to p \vee \neg p) \to \neg\neg(\neg p \vee \neg\neg p) \\
			\scriptstyle\equiv &\scriptstyle\neg\neg( (\neg\neg p \to p) \to p \vee \neg p) \to \top \\[.3em]
			\scriptstyle\equiv &\top
		\end{array}
		\end{equation*}
	\end{minipage}\hspace{1em}%
	\begin{minipage}[t]{0.48\textwidth}
		$\beta_5$ has property 1:
		\vspace{-.5em}
		\begin{equation*}
		\begin{array}{rl}
			\neg\neg \beta_5  \quad
			\scriptstyle\equiv &\scriptstyle\neg\neg( ((\neg\neg p \to p) \to p \vee \neg p) \vee (\neg\neg p \to p) ) \\
			\scriptstyle\equiv &\scriptstyle\neg\neg( \neg\neg ((\neg\neg p \to p) \to p \vee \neg p) \vee \neg\neg (\neg\neg p \to p) ) \\
			\scriptstyle\equiv &\scriptstyle\neg\neg( \neg\neg ((\neg\neg p \to p) \to p \vee \neg p) \vee \top ) \\
			\scriptstyle\equiv &\scriptstyle\neg\neg \top  \\[.3em]
			\scriptstyle\equiv &\top
		\end{array}
		\end{equation*}
	\end{minipage}

	\vspace{.3em}
	\noindent
	\begin{minipage}[t]{0.48\textwidth}
		$\alpha_5$ has property 2:
		for $\phi$ with property 1,
		\vspace{-.5em}
		\begin{equation*}
		\begin{array}{rl}
			\alpha_5(\phi)  \quad
			\scriptstyle\equiv &\scriptstyle((\neg\neg \phi \to \phi) \to \phi \vee \neg \phi) \to (\neg \phi \vee \neg\neg \phi) \\
			\scriptstyle\equiv &\scriptstyle((\neg\neg \phi \to \phi) \to \phi \vee \neg \phi) \to (\bot \vee \top) \\
			\scriptstyle\equiv &\scriptstyle((\neg\neg \phi \to \phi) \to \phi \vee \neg \phi) \to \top \\[.3em]
			\scriptstyle\equiv &\top
		\end{array}
		\end{equation*}
	\end{minipage}\hspace{1em}%
	\begin{minipage}[t]{0.48\textwidth}
		$\beta_5$ has property 2:
		for $\phi$ with property 1,
		\vspace{-.5em}
		\begin{equation*}
		\begin{array}{rl}
			\beta_5(\phi)  \quad
			\scriptstyle\equiv &\scriptstyle((\neg\neg \phi \to \phi) \to \phi \vee \neg \phi) \vee (\neg\neg \phi \to \phi) \\
			\scriptstyle\equiv &\scriptstyle((\top \to \phi) \to \phi \vee \bot) \vee (\neg\neg \phi \to \phi) \\
			\scriptstyle\equiv &\scriptstyle(\phi \to \phi) \vee (\neg\neg \phi \to \phi) \\
			\scriptstyle\equiv &\scriptstyle\top \vee (\neg\neg \phi \to \phi) \\[.3em]
			\scriptstyle\equiv &\top
		\end{array}
		\end{equation*}
	\end{minipage}

	\vspace{.3em}
	\noindent
	Moreover, we can show that, if $\alpha_n$ and $\beta_n$ have properties both properties, than this holds for $\alpha_{n+1}$ and $\beta_{n+1}$ too.

	\vspace{.3em}
	\noindent
	\begin{minipage}[t]{0.48\textwidth}
		$\alpha_{n+1}$ has property 1:
		\vspace{-.5em}
		\begin{equation*}
		\begin{array}{rl}
			\neg\neg \alpha_{n+1}  \quad
			\scriptstyle\equiv &\scriptstyle\neg\neg(\alpha_n \to \beta_n) \\
			\scriptstyle\equiv &\scriptstyle\neg\neg \alpha_n \to \neg\neg \beta_n \\
			\scriptstyle\equiv &\scriptstyle\neg\neg \alpha_n \to \top \\[.3em]
			\scriptstyle\equiv &\top
		\end{array}
		\end{equation*}
	\end{minipage}\hspace{1em}%
	\begin{minipage}[t]{0.48\textwidth}
		$\beta_{n+1}$ has property 1:
		\vspace{-.5em}
		\begin{equation*}
		\begin{array}{rl}
			\neg\neg \beta_{n+1}  \quad
			\scriptstyle\equiv &\scriptstyle\neg\neg( \alpha_n \vee \beta_n ) \\
			\scriptstyle\equiv &\scriptstyle\neg\neg( \neg\neg \alpha_n \vee \neg\neg \beta_n ) \\
			\scriptstyle\equiv &\scriptstyle\neg\neg( \top \vee \top ) \\[.3em]
			\scriptstyle\equiv &\top
		\end{array}
		\end{equation*}
	\end{minipage}

	\vspace{.3em}
	\noindent
	\begin{minipage}[t]{0.48\textwidth}
		$\alpha_{n+1}$ has property 2:
		for $\phi$ with property 1,
		\vspace{-.5em}
		\begin{equation*}
		\begin{array}{rl}
			\alpha_{n+1}(\phi)  \quad
			\scriptstyle\equiv &\scriptstyle\alpha_n(\phi) \to \beta_n(\phi) \\
			\scriptstyle\equiv &\scriptstyle\alpha_n(\phi) \to \top \\[.3em]
			\scriptstyle\equiv &\top
		\end{array}
		\end{equation*}
	\end{minipage}\hspace{1em}
	\begin{minipage}[t]{0.48\textwidth}
		$\beta_{n+1}$ has property 2:
		for $\phi$ with property 1,
		\vspace{-.5em}
		\begin{equation*}
		\begin{array}{rl}
			\beta_{n+1}(\phi)  \quad
			\scriptstyle\equiv &\scriptstyle\beta_n(\phi) \vee \alpha_n(\phi) \\
			\scriptstyle\equiv &\scriptstyle\top \vee \top \\[.3em]
			\scriptstyle\equiv &\top
		\end{array}
		\end{equation*}
	\end{minipage}

	\vspace{.3em}
	\noindent
	So, by induction all the formulas $\alpha_n, \beta_n$ with $n\geq 5$ have index at most $2$ and fixpoint $\top$.
	As for the remaining formulas, one can easily show their fix-points are as follows:
	\begin{equation*}
	\begin{array}{r@{\hspace{.3em}}ll @{\hspace{4em}\implies\hspace{4em}} r@{\hspace{.5em}}c@{\hspace{.5em}}l}
		\beta_1^2 &= (p)^2 &\equiv p   
			&(\beta_1^2)^0 &\equiv&(\beta_1^2)^1   \\
		\beta_2^2 &= (p \vee \neg p)^2 &\equiv p \vee \neg p   
			&(\beta_2^2)^1 &\equiv &(\beta_2^2)^3   \\
		\beta_3^2 &\equiv ( \neg p \vee \neg\neg p )^2 &\equiv \top   
			&(\beta_3^2)^2 &\equiv &(\beta_3^2)^4   \\
		\beta_4^2 &\equiv ( \neg\neg p \vee (\neg\neg p \to p) )^2 &\equiv \top   
			&(\beta_4^2)^2 &\equiv &(\beta_4^2)^4   \\
		\alpha_3^2 &\equiv ( \neg\neg p \to p )^2 &\equiv  \neg\neg p \to p   
			&(\alpha_3^2)^1 &\equiv &(\alpha_3^2)^2   \\
		\alpha_2^2 &= (\neg \neg p)^2 &\equiv \neg \neg p   
			&(\alpha_2^2)^1 &\equiv &(\alpha_2^2)^2   \\
		\alpha_1^3 &= (\neg p)^3 &\equiv \neg p   
			&(\alpha_1^3)^1 &\equiv &(\alpha_1^3)^3   \\
		\alpha_4^2 &\equiv ( (\neg\neg p \to p) \to p \vee \neg p )^2 &\equiv \top   
			&(\alpha_4^2)^2 &\equiv &(\alpha_4^2)^4   
	\end{array}
	\end{equation*}

\end{proof}

\end{document}